\documentclass[a4paper,10pt,twoside]{article}

\usepackage[utf8]{inputenc}
\usepackage{amsmath, amssymb, amsthm}
\usepackage{amstext, amsfonts, a4}
\usepackage{dsfont}
\usepackage{latexsym}
\usepackage{mathrsfs}
\usepackage[all,cmtip]{xy}
\usepackage{graphicx}
\usepackage{enumerate}
\usepackage{hyperref,color}
\usepackage{stmaryrd}
\usepackage[shortlabels]{enumitem}
\usepackage{bbm}

\usepackage{url}
\usepackage[pdftex,a4paper,left=3cm,right=3cm,top=3cm,bottom=3cm]{geometry}

\def\ra{\rightarrow}

\def\lra{\longrightarrow}

\def\lmapsto{\longmapsto}

 \def\sS{\mathscr{S}}

\def\bbA{\mathbb{A}}
\def\bbF{\mathbb{F}}\def\bbG{\mathbb{G}}

\def\bbP{\mathbb{P}}
\def\bbQ{\mathbb{Q}}\def\bbT{\mathbb{T}}

\def\cA{\mathcal{A}}
\def\cH{\mathcal{H}}

\def\cM{\mathcal{M}}\def\cO{\mathcal{O}}

\def\cX{\mathcal{X}}
\def\cZ{\mathcal{Z}}

\def\bfG{\mathbf{G}}

\def\bfT{\mathbf{T}}

\def\bfq{\mathbf{q}}

\def\whG{\widehat{G}}

\def\whT{\widehat{T}}

\DeclareMathOperator{\ASph}{ASph}

\DeclareMathOperator{\diag}{diag}

\DeclareMathOperator{\Ind}{Ind}
\DeclareMathOperator{\ind}{ind}

\DeclareMathOperator{\LP}{LP}
\DeclareMathOperator{\LzPz}{L_{\zeta}P_{\zeta}}

\DeclareMathOperator{\Mod}{Mod}

\DeclareMathOperator{\nonreg}{non-reg}

\DeclareMathOperator{\pr}{pr}

\DeclareMathOperator{\QCoh}{QCoh}

\DeclareMathOperator{\reg}{reg}

\DeclareMathOperator{\SingDiag}{SingDiag}

\DeclareMathOperator{\Spec}{Spec}

\DeclareMathOperator{\SP}{SP}

\DeclareMathOperator{\sss}{ss}

\DeclareMathOperator{\Sym}{Sym}

\DeclareMathOperator{\unr}{unr}

\newtheorem{counter}[subsection]{$\!\!$}

\newtheorem{counter*}[subsubsection]{$\!\!$}
\newenvironment{Def*}{\begin{counter*} {\bf Definition.}}{\end{counter*}}
\newenvironment{Not*}{\begin{counter*} \rm {\bf Notation.}}{\end{counter*}}
\newenvironment{Notss*}{\begin{counter*} \rm {\bf Notations.}}{\end{counter*}}
\newenvironment{DefNot*}{\begin{counter*} \rm {\bf Definition-Notation.}}{\end{counter*}}
\newenvironment{Nots*}{\begin{counter*} \rm {\bf Notations.}}{\end{counter*}}
\newenvironment{Prop*}{\begin{counter*} {\bf Proposition.}}{\end{counter*}}
\newenvironment{Lem*}{\begin{counter*} {\bf Lemma.}}{\end{counter*}}
\newenvironment{Cor*}{\begin{counter*} {\bf Corollary.}}{\end{counter*}}
\newenvironment{Th*}{\begin{counter*} {\bf Theorem.}}{\end{counter*}}
\newenvironment{Rem*}{\begin{counter*} \rm {\bf Remark.}}{\end{counter*}}
\newenvironment{Ex*}{\begin{counter*} \rm {\bf Example.}}{\end{counter*}}
\newenvironment{Exs*}{\begin{counter*} \rm {\bf Examples.}}{\end{counter*}}
\newenvironment{Pt*}{\begin{counter*} \rm}{\end{counter*}}
\newenvironment{Q*}{\begin{counter*} \rm {\bf Question.}}{\end{counter*}}

\newcommand{\iso}{\stackrel{\sim}{\longrightarrow}}
\newcommand{\EG}{X_{\zeta}}

\title{\textbf{\huge{A semisimple mod $p$ Langlands correspondence in families for $GL_2(\bbQ_p)$}}}

\author{Cédric PEPIN and Tobias SCHMIDT}
\date{\today}

\begin{document}

\maketitle

\begin{abstract}
This is the sequel to \cite{PS}. Let $F$ be any local field with residue characteristic $p>0$, and $\cH^{(1)}_{\overline{\bbF}_p}$ be the mod $p$ pro-$p$-Iwahori Hecke algebra of $\mathbf{GL_2}(F)$. In \cite{PS} we have constructed a parametrization of the $\cH^{(1)}_{\overline{\bbF}_p}$-modules by certain $\widehat{\mathbf{GL_2}}(\overline{\bbF}_p)$-Satake parameters, together with an antispherical family of $\cH^{(1)}_{\overline{\bbF}_p}$-modules. Here we let $F=\bbQ_p$ (and $p\geq 5$) and construct a morphism from $\widehat{\mathbf{GL_2}}(\overline{\bbF}_p)$-Satake parameters to $\widehat{\mathbf{GL_2}}(\overline{\bbF}_p)$-Langlands parameters. 
As a result, we get a version in families of Breuil's semisimple mod $p$ Langlands correspondence for $\mathbf{GL_2}(\bbQ_p)$ and of Pa\v{s}k\={u}nas' parametrization of blocks of the category of mod $p$ locally admissible smooth representations of $\mathbf{GL_2}(\bbQ_p)$ having a central character. The formulation of these results is possible thanks to the Emerton-Gee moduli space of semisimple $\widehat{\mathbf{GL_2}}(\overline{\bbF}_p)$-representations of the Galois group ${\rm Gal}(\overline{\bbQ}_p/ \bbQ_p)$. 
\end{abstract}

\tableofcontents

\section*{Introduction}
Let $F$ be a local field with ring of integers $o_F$ and residue field $\bbF_q$. We let $\bfG$ be the algebraic group $\mathbf{GL_2}$ over $F$ with diagonal torus $\bfT\subset\bfG$. Set $G:=\bfG(F)$. Let $\cH^{(1)}_{\overline{\bbF}_q}$ be the pro-$p$-Iwahori Hecke algebra of $G$, with coefficients in an algebraic closure $ \overline{\bbF}_q$ of $\bbF_q$. 
Let $\widehat{\bfG}$ be the Langlands dual group of $\bfG$ over $\overline{\bbF}_q$, with maximal torus $\widehat{\bfT}$. In this sequel to \cite{PS}, we continue to work at $\bfq=q=0$.  
That is, we consider the special fibre at $\bfq=0$ of the Vinberg fibration $V_{\mathbf{\whT}}\stackrel{\bfq}{\rightarrow}\bbA^1$ associated to $\widehat{\bfT}$ followed by base change to $\overline{\bbF}_q$. This yields
the ${\overline{\bbF}_q}$-semigroup scheme

$$
V_{\mathbf{\whT},0} := \SingDiag_{2\times2}\times_{\overline{\bbF}_q} \bbG_m,
$$
where $\SingDiag_{2\times2}$ represents the semigroup of singular diagonal $2\times 2$-matrices over $\overline{\bbF}_q$, cf. \cite[7.1]{PS}. Let $\bbT^{\vee}$ be the finite abelian dual group of $\bbT=\bfT(\bbF_q)$ and consider the extended semigroup
$$V^{(1)}_{\mathbf{\whT},0}:=\bbT^{\vee}\times V_{\mathbf{\whT},0}.$$ 
It has a natural $W_0$-action. In \cite[7.2.2]{PS} we established the mod $p$ pro-$p$-Iwahori Satake isomorphism
$$
\xymatrix{
\sS^{(1)}_{\overline{\bbF}_q}: Z(\cH^{(1)}_{\overline{\bbF}_q})\ar[r]^<<<<<{\sim} & \cO(V^{(1)}_{\mathbf{\whT},0}/W_0)
}
$$
identifying the center $Z(\cH^{(1)}_{\overline{\bbF}_q})\subset \cH^{(1)}_{\overline{\bbF}_q}$ with the ring of regular functions 
on the quotient $V^{(1)}_{\mathbf{\whT},0}/W_0$. The resulting Satake equivalence $S$ identifies the category of $Z(\cH^{(1)}_{\overline{\bbF}_q})$-modules with the category of $\widehat{\bfG}$-Satake parameters, i.e. the category of quasi-coherent sheaves on $V^{(1)}_{\mathbf{\whT},0}/W_0$, cf. \cite[7.3.2]{PS}.

We also constructed the mod $p$ antispherical module $\cM^{(1)}_{\overline{\bbF}_q}$, cf. \cite[7.4.1]{PS}. This is a distinguished $\cH^{(1)}_{\overline{\bbF}_q}$-action on the maximal 
commutative subring $\cA^{(1)}_{\overline{\bbF}_q}$ of $\cH^{(1)}_{\overline{\bbF}_q}$. The sheaf
$S(\cM^{(1)}_{\overline{\bbF}_q})$,
when specialized at closed points of $V^{(1)}_{\mathbf{\whT},0}/W_0$, gives rise to 
a dual parametrization of {\it all} irreducible $\cH^{(1)}_{\overline{\bbF}_q}$-modules in terms of $\widehat{\bfG}$-Satake parameters \cite[7.4.9/7.4.15]{PS}.

\vskip5pt

In this sequel to \cite{PS} we construct, in the case $F=\bbQ_p$ and $p\geq 5$, a morphism $L$ from the space of $\widehat{\bfG}$-Satake parameters to the space of $\widehat{\bfG}$-Langlands parameters, and prove that the push-forward $L_*S(\cM^{(1)}_{\overline{\bbF}_q})$ interpolates the semisimple mod $p$ local Langlands correspondence $\rho\mapsto \pi(\rho)$ for the group $G$. 

\vskip5pt 

To be more precise, let $\zeta: Z(G)\rightarrow {\overline{\bbF}^{\times}_q}$ be a central character of $G$. There is a natural fibration $\theta: V^{(1)}_{\mathbf{\whT},0}/W_0\rightarrow Z(G)^\vee$
where $Z(G)^{\vee}$ is the group scheme of characters of $Z(G)$, and we put
$$(V^{(1)}_{\mathbf{\whT},0}/W_0)_{\zeta}:= \theta^{-1}(\zeta).$$ 

We let from now on $F=\bbQ_p$ with $p\geq 5$. As a space of $\widehat{\bfG}$-Langlands parameters, we may then consider the Emerton-Gee moduli curve $X_{\zeta}$, cf. \cite{Em19}, parametrizing (isomorphism classes of) two-dimensional semisimple continuous Galois representations over $\overline{\bbF}_p$ with determinant $\omega\zeta$:
$$
\EG(\overline{\bbF}_p) \cong \big\{\text{semisimple continuous $\rho : {\rm Gal}(\overline{\bbQ}_p/ \bbQ_p)\rightarrow \mathbf{\whG}(\overline{\bbF}_p)$\;
with $\det \rho = \omega \zeta$} \big\}/ \sim.
$$
Here $\omega$ is the mod $p$ cyclotomic character. 
The curve $X_{\zeta}$ is expected to be the underlying scheme of a ringed moduli space
for the stack of \'etale $(\varphi,\Gamma)$-modules $\cX_2^{\det = \omega\zeta}$ appearing in \cite{EG19} (see also \cite{CEGS19}). At the moment, it is unclear how to define a replacement for $X_{\zeta}$ when $F/\bbQ_p$ is a non trivial finite extension, and this is the reason why we have to restrict to the case $F=\bbQ_p$ (and $p\geq 5$) in our construction of the morphism $L$. Our main result is the following (cf. Theorem \ref{LLfamily}).

\vskip5pt 

{\bf Theorem.}{\it  \; Suppose $F=\bbQ_p$ with $p\geq 5$. There exists a morphism of $\overline{\bbF}_p$-schemes  
$$ L_\zeta: (V^{(1)}_{\mathbf{\whT},0}/W_0)_{\zeta}\longrightarrow \EG$$ 
such that the quasi-coherent $\cO_{X_{\zeta}}$-module 
$$
L_{\zeta*}S(\cM_{\overline{\bbF}_p}^{(1)})|_{(V^{(1)}_{\mathbf{\whT},0}/W_0)_{\zeta}},
$$ 
equal to the push-forward along $L_{\zeta}$ of the restriction to $(V^{(1)}_{\mathbf{\whT},0}/W_0)_{\zeta}\subset V^{(1)}_{\mathbf{\whT},0}/W_0$ of the Satake parameter $S(\cM_{\overline{\bbF}_p}^{(1)})$, interpolates the $I^{(1)}$-invariants of the semisimple mod $p$ Langlands correspondence
$$
\left.
\begin{array}{lllll}
X_{\zeta}(\overline{\bbF}_p) & \lra & \Mod_{\zeta}^{\rm l adm}(\overline{\bbF}_p[G]) & \lra & \Mod(\cH^{(1)}_{\overline{\bbF}_p}) \\
x & \lmapsto & \pi(\rho_x) & \lmapsto & \pi(\rho_x)^{I^{(1)}},
\end{array}
\right.
$$
in the sense: for all $x\in X_{\zeta}(\overline{\bbF}_p)$, one has an isomorphism of $\cH^{(1)}_{\overline{\bbF}_p}$-modules 

$$
\Big(\big(L_{\zeta*}S(\cM_{\overline{\bbF}_p}^{(1)})|_{(V^{(1)}_{\mathbf{\whT},0}/W_0)_{\zeta}}\big)\otimes_{\cO_{X_{\zeta}}}k(x)\Big)^{\sss}=\Big(\cM_{\overline{\bbF}_p}^{(1)}\otimes_{Z(\cH^{(1)}_{\overline{\bbF}_p})}(\sS_{\overline{\bbF}_p}^{(1)})^{-1}(\cO_{L_{\zeta}^{-1}(x)})\Big)^{\sss}\cong \pi(\rho_x)^{I^{(1)}}.
$$
}

\vskip5pt

Here, $\Mod_{\zeta}^{\rm l adm}(\overline{\bbF}_p[G]) $ denotes the category of locally admissible smooth $G$-representations over $\overline{\bbF}_p$ with central character $\zeta$. The group $I^{(1)}\subset G$ is the standard pro-$p$ Iwahori subgroup and $(\cdot)^{I^{(1)}}$ denotes the functor of $I^{(1)}$-invariants.

\vskip5pt 

As a byproduct of our constructions, we also obtain a version in families of Pa\v{s}k\={u}nas' parametrization of the blocks of the category 
$\Mod_{\zeta}^{\rm l adm}(\overline{\bbF}_p[G])$, cf. \cite{Pas13}. See \ref{paramfamily} for the precise statement.

\section{The theory at $\bfq=q=0$: Semisimple Langlands correspondence}

We keep the notation from the introduction. 
In particular, $F$ denotes a local field with ring of integers $o_F$ and residue field $\bbF_q$
(we switch to $F=\bbQ_p$ starting from \ref{LPfixdet}). We also let $k:=\overline{\bbF}_q$.


\subsection{Mod $p$ Satake parameters with fixed central character}\label{Sat_par_with_cc}

\begin{Pt*} \label{Fq_action} Let $\omega : \bbF_q^{\times}\ra k^\times$ be induced by the inclusion $\bbF_q \subset k$. 
Then $(\bbF_q^{\times})^{\vee}=\langle \omega \rangle$ is a cyclic group of order $q-1$. 
An element $\omega^r$ defines a non-regular character of $\bbT$: 
$$\omega^r(t_1,t_2):=\omega^{r}(t_1)\omega^{r}(t_2)$$ 
for all $(t_1,t_2)\in\bbT=\bbF_q^{\times}\times\bbF_q^{\times}$. Composing with multiplication in $\bbT^{\vee}$, we get an action of 
$(\bbF_q^{\times})^{\vee}$ on  $\bbT^{\vee}$,
which factors on the quotient set $\bbT^{\vee}/W_0$:
$$  \bbT^{\vee}/W_0 \times (\bbF_q^{\times})^{\vee} \longrightarrow  \bbT^{\vee}/W_0, \;(\gamma, \omega^r) \mapsto \gamma\omega^r.$$
If $\gamma\in \bbT^{\vee}/W_0$ is regular (non-regular), then $\gamma\omega^r$
is regular (non-regular). 
\end{Pt*}

\begin{Pt*} \label{orbits}
Restricting characters of $\bbT$ to the subgroup $\bbF_q^{\times}\simeq \{ \diag(a,a) : a\in \bbF_q^{\times}\} $ induces a homomorphism $\bbT^{\vee}\ra   (\bbF_q^{\times})^{\vee}$ which factors into a restriction map  
$$ \bbT^{\vee}/W_0 \ra   (\bbF_q^{\times})^{\vee}, \; \gamma\mapsto \gamma |_{\bbF_q^{\times}}.$$
The relation to the $(\bbF_q^{\times})^{\vee} $-action on the source $\bbT^{\vee}/W_0$ is given by the formula
$$(\gamma\omega^r)|_{ \bbF_q^{\times}} = \gamma|_{ \bbF_q^{\times}} \;\omega^{2r}.$$ 

We describe the fibers of the restriction map $\gamma\mapsto \gamma |_{\bbF_q^{\times}}$.

Let $(\cdot)|_{ \bbF_q^{\times}}^{-1}(\omega^{2r})$ be the fibre at a square element
$\omega^{2r}$. By the above formula, the action of $\omega^{-r}$ on $\bbT^{\vee}/W_0$ induces a bijection with the fibre $(\cdot)|_{\bbF_q^{\times}}^{-1}(1).$
The fibre 
$$(\cdot)|_{\bbF_q^{\times}}^{-1}(1)=\{ 1 \otimes 1 \} \coprod \{ \omega \otimes \omega^{-1}, \omega^2 \otimes \omega^{-2},...,  
\omega^{\frac{q-3}{2}} \otimes \omega^{-{\frac{q-3}{2}}} \} \coprod \{ \omega^{\frac{q-1}{2}} \otimes \omega^{-{\frac{q-1}{2}}} \}$$
has cardinality $\frac{q+1}{2}$ and, in the above list, we have chosen a representative in $\bbT^{\vee}$ for each element in the fibre. 
The $\frac{q-3}{2}$ elements in the middle of this list, i.e. the $W_0$-orbits 
represented by the characters $\omega^r \otimes \omega^{-r}$ for $r=1,...,\frac{q-3}{2}$, are all regular $W_0$-orbits. The two orbits at the two ends of the list
are non-regular orbits (note that $\frac{q-1}{2}\equiv -\frac{q-1}{2} \mod (q-1)$). Since the action of $\omega^{-r}$ preserves regular (non-regular) orbits, any fibre at a square element (there are 
$\frac{q-1}{2}$ such fibres) has the same structure.

On the other hand, let $(\cdot)|_{\bbF_q^{\times}}^{-1}(\omega^{2r-1})$ be the fibre at a non-square element $\omega^{2r-1}$. The action of 
$\omega^{-r}$ induces a bijection 
with the fibre $(\cdot)|_{\bbF_q^{\times}}^{-1}(\omega^{-1})$.  
The fibre 
$$(\cdot)|_{\bbF_q^{\times}}^{-1}(\omega^{-1})= \{ 1 \otimes \omega^{-1}, \omega \otimes \omega^{-2},...,  
\omega^{\frac{q-1}{2}-1} \otimes \omega^{-{\frac{q-1}{2}}} \} $$
has cardinality $\frac{q-1}{2}$ and we have chosen a representative in $\bbT^{\vee}$ for each element in the fibre. 
All elements of the fibre are regular $W_0$-orbits. Since the action of $\omega^{-r}$ preserves regular (non-regular) orbits, any fibre at a non-square element (there are $\frac{q-1}{2}$ such fibres) has the same structure.

Note that $\frac{q-1}{2}( \frac{q+1}{2} +  \frac{q-1}{2} ) = \frac{q^2-q}{2}$ is the cardinality of the set  $\bbT^{\vee}/W_0$.
\end{Pt*}

\begin{Pt*}\label{ZGveeaction}
Recall the commutative $k$-semigroup scheme
$$
V^{(1)}_{\mathbf{\whT},0}=\bbT^{\vee}\times V_{\mathbf{\whT},0}=\bbT^{\vee}\times  \SingDiag_ {2\times 2}\times\bbG_m 
$$ 
together with its $W_0$-action, cf. \cite[6.2.15]{PS}: the natural action of $W_0$ on 
the factors $\bbT^{\vee}$ and  $\SingDiag_ {2\times 2}$ and the trivial one on $\bbG_m$. 
There is a commuting action of the $k$-group scheme
$$\cZ^\vee:= (\bbF_q^{\times})^{\vee}\times \bbG_m$$ 
on $V^{(1)}_{\mathbf{\whT},0}$:  
the (constant finite diagonalizable) group $(\bbF_q^{\times})^{\vee}$ acts only on the factor $\bbT^{\vee}$ and in the way described in \ref{Fq_action}; an element $z_0\in\bbG_m$ acts trivially on $\bbT^{\vee}$, by multiplication with the diagonal matrix $\diag(z_0,z_0)$ on $\SingDiag_ {2\times 2}$ and by multiplication with the square $z_0^2$ on $\bbG_m$. Therefore the quotient 
$V^{(1)}_{\mathbf{\whT},0}/W_0$ inherits a $\cZ^\vee$-action. Now, according to \cite[7.4.7]{PS}, one has the decomposition
$$
V^{(1)}_{\mathbf{\whT},0}/W_0=\coprod_{\gamma\in(\bbT^{\vee}/W_0)_{\reg}} V_{\mathbf{\whT},0}\coprod_{\gamma\in(\bbT^{\vee}/W_0)_{\nonreg} } V_{\mathbf{\whT},0}/W_0.
$$
Then the $(\bbF_q^{\times})^{\vee}$-action is by permutations on the index set $\bbT^{\vee}/W_0$, i.e. on the 
set of connected components of $V^{(1)}_{\mathbf{\whT},0}/W_0$; as observed above, it preserves the subsets of regular and non-regular components. The $\bbG_m$-action on $V^{(1)}_{\mathbf{\whT},0}/W_0$ preserves each connected component.

\end{Pt*} 
\begin{Pt*}\label{twistHecke} Recall from \cite[7.4.7]{PS} the antispherical map 
 $$
\xymatrix{
\ASph:  (V^{(1)}_{\mathbf{\whT},0}/W_0)(k)\ar[r] & \{\textrm{left $\cH^{(1)}_{\overline{\bbF}_q}$-modules}\}/\sim.
}
$$
The modules in the image of this map are standard modules of length $1$ or $2$, cf.  \cite[7.4.9]{PS} and  \cite[7.4.15]{PS}.

Let $(\omega^r,z_0)\in \cZ^\vee(k)$. Then recall that the standard $\cH^{(1)}_{\overline{\bbF}_q}$-modules and their simple constituents 
may be `twisted by the character $(\omega^r,z_0)$' : in the regular case, the actions of $X,Y,U^2$ get multiplied by $z_0,z_0,z_0^{2}$ respectively and the component $\gamma$ gets multiplied by $\omega^r$, cf. \cite[2.4]{V04}; in the non-regular case, the action of $U$ gets multiplied by $z_0$, the action of $S$ remains unchanged and the component $\gamma$ gets multiplied by $\omega^r$, cf. \cite[1.6]{V04}.
This gives an action of the group of $k$-points of $\cZ^\vee$ on the standard $\cH^{(1)}_{\overline{\bbF}_q}$-modules and their simple constituents.

\end{Pt*}
\begin{Lem*} \label{Asphequiv} The map $\ASph$ is $\cZ^\vee(k)$-equivariant. 
\end{Lem*}
\begin{proof} Let $(\omega^r,z_0)\in \cZ^\vee(k)$.
Let $v\in  (V^{(1)}_{\mathbf{\whT},0}/W_0)(k)$ and let its connected component be indexed by $\gamma\in\bbT^{\vee} /W_0$. 
Suppose that $\gamma$ is regular, choose an ordering $\gamma=(\chi,\chi^s)$ on the set $\gamma$ and standard coordinates. Then
$\ASph(v)=\ASph^{\gamma}(v)$ is a simple two-dimensional standard $\cH^{\gamma}_{\overline{\bbF}_p}$-module, cf.   \cite[7.4.9]{PS}, i.e. of the form
$M(x,y,z_2,\chi)$ \cite[3.2]{V04}. Then $$\ASph(v.(\omega^r,z_0))\simeq M(z_0x,z_0y,z_0^2z_2,\chi. \omega^r)\simeq\ASph(v).(\omega^r,z_0).$$
Suppose that $\gamma=\{\chi\}$ is non-regular and choose Steinberg coordinates. (a) If $v\in D(2)_{\gamma}(k)$, then 
$\ASph(v)=\ASph^{\gamma}(2)(v)$ is 
a simple two-dimensional $\cH^{\gamma}_{\overline{\bbF}_p}$-module, cf.  \cite[7.4.15]{PS}, i.e. of the form 
$M(z_1,z_2,\chi)$ \cite[3.2]{V04}. 
Then $$\ASph(v.(\omega^r,z_0))\simeq M(z_0z_1,z_0^2z_2,\chi. \omega^r)\simeq \ASph(v).(\omega^r,z_0).$$
(b) If  $v\in D(1)_{\gamma}(k)$, then the semisimplified module $\ASph(v)^{\sss}$ is the direct sum of the two characters in the antispherical pair 
$\ASph^{\gamma}(1)(v)=\{(0,z_1),(-1,-z_1)\}$ where $z_2=z^2_1$. Similarly $\ASph(v.(\omega^r,z_0))^{\sss}$ is the direct sum of the 
characters $\{(0,z_0z_1),(-1,-z_0z_1)\}$ in the component $\gamma.\omega^r$, and hence is isomorphic to $\ASph(v)^{\sss}.(\omega^r,z_0)$.
\end{proof}

\begin{Pt*} \label{projections}
The two canonical projections from $V^{(1)}_{\mathbf{\whT},0}$ to $\bbT^{\vee}$ and 
$\bbG_m$ respectively induce two projection morphisms 

$$
\xymatrix{
&V^{(1)}_{\mathbf{\whT},0}/W_0 \ar[dl]_{\pr_{\bbT^{\vee}/W_0}} \ar[dr]^{\pr_{\bbG_m}} &  \\
\bbT^{\vee}/W_0  & &\bbG_m .
}
$$
Then we may compose the map $\pr_{\bbT^{\vee}/W_0}$ with the restriction map $(\cdot) |_{\bbF^{\times}_q} :\bbT^{\vee}/W_0\ra(\bbF_q^{\times})^{\vee}$, set
$$\theta := \big((\cdot) |_{\bbF^{\times}_q} \circ \pr_{\bbT^{\vee}/W_0}\big) \times {\pr_{\bbG_m}} $$
and view $V^{(1)}_{\mathbf{\whT},0}/W_0$ as fibered over the space $\cZ^\vee$:
$$
\xymatrix{
&V^{(1)}_{\mathbf{\whT},0}/W_0 \ar[d]^\theta &  \\
&\cZ^\vee. &
}
$$
The relation to the $\cZ^\vee$-action on the source $V^{(1)}_{\mathbf{\whT},0}/W_0$ is given by the formula
$$\theta (x.(\omega^r,z_0) )=\theta(x)(\omega^{2r},z^2_0)=\theta(x)(\omega^{r},z_0)^2$$ 
for $x\in V^{(1)}_{\mathbf{\whT},0}/W_0$ and $(\omega^r,z_0)\in\cZ^{\vee}$. This formula 
follows from the formula in \ref{orbits} and the definition of the $\bbG_m$-action in \ref{ZGveeaction}.
\end{Pt*}

\begin{Def*} \label{def_space_Sat_cc}
Let $\zeta\in \cZ^\vee$. The \emph{space of mod $p$ Satake parameters with central character $\zeta$} 
is the $k$-scheme $$(V^{(1)}_{\mathbf{\whT},0}/W_0)_{\zeta}:= \theta^{-1}(\zeta).$$
\end{Def*}

\begin{Pt*}\label{connectedcomp}
Let $\zeta=(\zeta |_ {\bbF^{\times}_q}, z_2)\in \cZ^\vee(k)=(\bbF^{\times}_q)^{\vee}\times k^{\times}$. Denote by $(V^{(1)}_{\mathbf{\whT},0}/W_0)_{z_2}$ the fibre of $\pr_{\bbG_m}$ at $z_2\in k^{\times}$. Then by  \cite[7.4.7]{PS} we have
$$
(V^{(1)}_{\mathbf{\whT},0}/W_0)_{\zeta}=\coprod_{\gamma\in(\bbT^{\vee}/W_0)_{\reg},  \gamma |_ {\bbF^{\times}_q}= \zeta |_ {\bbF^{\times}_q} } V_{\mathbf{\whT},0,z_2}\coprod_{\gamma\in(\bbT^{\vee}/W_0)_{\nonreg}, \gamma |_ {\bbF^{\times}_q}= \zeta |_ {\bbF^{\times}_q} } V_{\mathbf{\whT},0,z_2}/W_0.
$$
Recall that the choice of standard coordinates $x,y$ identifies
$$
V_{\mathbf{\whT},0,z_2}\simeq \bbA^1\cup_0 \bbA^1
$$ 
with two affine lines over $k$, intersecting at the origin, cf.  \cite[7.4.8]{PS}.
On the other hand, the choice of the Steinberg coordinate $z_1$ identifies
$$
V_{\mathbf{\whT},0,z_2}/W_0\simeq \bbA^1$$ 
with a single affine line over $k$, cf.  \cite[7.4.10]{PS}.
\end{Pt*}

\begin{Lem*}\label{twistSatake} Let $\zeta, \eta\in \cZ^{\vee}$. The action of $\eta$ on 
$V^{(1)}_{\mathbf{\whT},0}/W_0$ induces an isomorphism of $k$-schemes $(V^{(1)}_{\mathbf{\whT},0}/W_0)_{\zeta} \simeq (V^{(1)}_{\mathbf{\whT},0}/W_0)_{\zeta\eta^2}$.
\end{Lem*} 

\begin{proof} 
Follows from the last formula in \ref{projections}. 
\end{proof}

\subsection{Mod $p$ Langlands parameters with fixed determinant for $F=\bbQ_p$} \label{LPfixdet}

\begin{Not*}
In this section, we let $F=\bbQ_p$ with $p\geq 5$. 
We fix an algebraic closure $\overline{\bbQ}_p$ and let ${\rm Gal}(\overline{\bbQ}_p/ \bbQ_p)$ be the
absolute Galois group. We normalize local class field theory $\bbQ_p^{\times}\rightarrow {\rm Gal}(\overline{\bbQ}_p/ \bbQ_p)^{\rm ab}$ by sending 
$p$ to a geometric Frobenius. In this way, we identify the $k$-valued smooth characters of ${\rm Gal}(\overline{\bbQ}_p/ \bbQ_p)$ and of $\bbQ_p^{\times}$.  Finally, $\omega: \bbQ_p^{\times} \ra k^\times$ denotes the extension of 
the character $\omega: \bbF_p^{\times}\ra k^{\times}$ to $\bbQ_p^{\times}$ satisfying $\omega(p)=1$, and $\unr(x): \bbQ_p^{\times} \ra k^\times$ denotes the character trivial on $\bbF_p^{\times}$ and sending $p$ to $x$. 

\end{Not*}
\begin{Pt*} \label{twistEG}
Let $\zeta:  \bbQ_p^{\times} \ra k^\times$ be a character. Recall from \cite{Em19} that the {\it Emerton-Gee moduli curve with character $\zeta$} is a certain projective curve $X_{\zeta}$ over $k$ whose points parametrize (isomorphism classes of) two-dimensional semisimple continuous Galois representations over $k$ with determinant $\omega\zeta$:
$$
\EG(k) \cong \big\{\text{semisimple continuous $\rho : {\rm Gal}(\overline{\bbQ}_p/ \bbQ_p)\rightarrow \mathbf{\whG}(k)$\;
with $\det \rho = \omega \zeta$} \big\}/ \sim.
$$
The curve $X_{\zeta}$ is a chain of projective lines over $k$ of length $\frac{p\pm 1}{2}$, whose irreducible components intersect at ordinary double points.
The sign $\pm 1$ is equal to $- \zeta (-1)$. We refer to $\zeta$ in the case $- \zeta (-1)=-1$ resp. $- \zeta (-1)=+1$ as an \emph{even character} resp. \emph{odd character}. 
There is a finite set of closed points $\EG^{\rm irred}\subset \EG$ which correspond to the classes of irreducible representations. 
Its open complement $\EG^{\rm red}=\EG\setminus  \EG^{\rm irred}$ parametrizes the reducible representations (i.e. direct sums of characters). Let  $\eta:   {\rm Gal}(\overline{\bbQ}_p/\bbQ_p) \ra k^\times$ be a character. Since $\det (\rho\otimes\eta) = (\det \rho)\eta^2$, twisting representations with $\eta$ induces an isomorphism 
$$(\cdot)\otimes\eta: \EG\iso X_{\zeta \eta^2}.$$
Hence one is reduced to consider only two `basic' cases: the even case where $\zeta(p)=1$ and $\zeta |_{\bbF_p^{\times}} = 1$ and the odd case where $\zeta(p)=1$ and $\zeta |_{\bbF_p^{\times}} = \omega^{-1}.$ Indeed, if $\zeta |_{\bbF_p^{\times}} = \omega^{r}$ for some even $r$, then choosing $\eta$ with $\eta(p)^2=\zeta(p)^{-1}$ and $\eta |_{\bbF_p^{\times}} = \omega^{-\frac{r}{2}}$,
one finds that $(\zeta\eta^2)(p)=1$ and $(\zeta\eta^2)  |_{\bbF_p^{\times}} = 1$; if $\zeta |_{\bbF_p^{\times}} = \omega^r$ for some odd $r$, then choosing $\eta$ with $\eta(p)^2=\zeta(p)^{-1}$ and $\eta |_{\bbF_p^{\times}} = \omega^{-\frac{r+1}{2}}$,
one finds that $(\zeta\eta^2)(p)=1$ and $(\zeta\eta^2)  |_{\bbF_p^{\times}} = \omega^{-1}$. 
\end{Pt*}

\begin{Pt*}  \label{structureEGeven} We make explicit some structure elements of $\EG$ in the even case $\zeta(p)=1$ and $\zeta |_{\bbF_p^{\times}} = 1$. 
Every irreducible component of $\EG$ is isomorphic to $\bbP^1$ and there are $\frac{p-1}{2}$ components. They are labelled by pairs of Serre weights of the following form: 
$$
\begin{array}{ccc}
\Sym^0  &  | &  \Sym^{p-3}\otimes\det  \\
 \Sym^2\otimes\det^{-1} & |  &  \Sym^{p-5}\otimes\det^2  \\
 \Sym^4\otimes\det^{-2} & |  &   \Sym^{p-7}\otimes\det^3 \\
 \vdots &  \vdots & \vdots   \\
 \Sym^{p-3}\otimes\det^{\frac{p+1}{2}}  & | &   \Sym^{0}\otimes\det^{\frac{p-1}{2}}.
\end{array}
$$
The component with label $"\Sym^0 |  \;   \Sym^{p-3}\otimes\det"$ intersects the next component at the point of $\EG^{\rm irred}$ 
parametrizing the irreducible Galois representation whose associated Serre weights are $\{ \Sym^2\otimes\det^{-1}, \Sym^{p-3}\otimes\det \}$. The component with label  $"\Sym^2\otimes\det^{-1} \; | \;  \Sym^{p-5}\otimes\det^2"$ intersects the next component at the point of 
$\EG^{\rm irred}$ parametrizing the irreducible Galois representation whose associated Serre weights are $\{ \Sym^4\otimes\det^{-2}, \Sym^{p-5}\otimes\det^2\}$. Continuing in this way, one finds $\frac{p-3}{2}$ points of $\EG^{\rm irred}$, which 
correspond to the $\frac{p-3}{2}$ double points of the chain $\EG$. There are two more points in $\EG^{\rm irred}$: they are smooth points, each one lies on one of the two `exterior' components and corresponds there to the irreducible Galois representation whose associated Serre weights are $\{ \Sym^0,\Sym^{p-1} \} $ and  $\{ \Sym^0\otimes\det^{\frac{p-1}{2}}, \Sym^{p-1}\otimes \det^{\frac{p-1}{2}} \} $ respectively. So $\EG^{\rm irred}$ has cardinality $\frac{p+1}{2}$.
 Suppose we are on one of the two exterior components $\bbP^1$. There is a canonical affine coordinate $z_1$ on the open complement of the double point, identifying this open complement with $\bbA^1$. We call the four points where $z_1=\pm 1$ {\it the four exceptional points} of $X_{\zeta}$.
\end{Pt*}

\begin{Pt*} \label{structureEGodd} We make explicit some structure elements of $\EG$ in the odd case $\zeta(p)=1$ and $\zeta |_{\bbF_p^{\times}} = \omega^{-1}$. Every irreducible component of $\EG$ is isomorphic to $\bbP^1$ and there are $\frac{p+1}{2}$ components. They are labelled by pairs of Serre weights of the following form: 
$$
\begin{array}{ccc}
\Sym^{p-2}  &  | &  "\Sym^{-1}"  \\
 \Sym^{p-4}\otimes\det & |  &  \Sym^1 \otimes\det^{-1}  \\
 \Sym^{p-6}\otimes\det^{2} & |  &   \Sym^{3}\otimes\det^{-2} \\
 \vdots &  \vdots & \vdots   \\
 \Sym^1\otimes\det^{\frac{p-3}{2}}  & | &   \Sym^{p-4}\otimes\det^{\frac{p+1}{2}} \\
 " \Sym^{-1}\otimes\det^{\frac{p-1}{2}}"  & | &   \Sym^{p-2}\otimes\det^{\frac{p-1}{2}}.
\end{array}
$$
The component with label $"\Sym^{p-2}  \;  | \;  "\Sym^{-1}" "$ intersects the next component at the point of $\EG^{\rm irred}$ 
parametrizing the irreducible Galois representation whose associated Serre weights are  $\{ \Sym^1\otimes\det^{-1}, \Sym^{p-2} \}$. The component with label  $" \Sym^{p-4}\otimes\det \; | \;   \Sym^1 \otimes\det^{-1}"$ intersects the next component at the point of $\EG^{\rm irred}$ parametrizing the irreducible Galois representation whose associated Serre weights are $\{ \Sym^3\otimes\det^{-2}, \Sym^{p-4}\otimes\det\}$. Continuing in this way, one finds $\frac{p-1}{2}$ points of $\EG^{\rm irred}$, which 
correspond to the $\frac{p-1}{2}$ double points of the chain $\EG$. There are no more points in $\EG^{\rm irred}$ and $\EG^{\rm irred}$ has cardinality $\frac{p-1}{2}$. Suppose we are on one of the two exterior components $\bbP^1$. There is a canonical affine coordinate $t$ on the open complement of the double point, identifying this open complement with $\bbA^1$. We call the four points where $t=\pm 2$ {\it the four exceptional points} of $X_{\zeta}$.
\footnote{The Galois representations living on the two exterior components in the odd case are {\it unramified} (up to twist), i.e. of type 
$
\rho=
\left (\begin{array}{cc}
\unr(x)& 0\\
0 & \unr(x^{-1})
\end{array} \right)\otimes\eta
$
and $t$ equals the `trace of Frobenius' $x+x^{-1}$. Hence $t=\pm 2$ if and only if $x=\pm 1$. }
\end{Pt*}

\begin{Def*}\label{modpLP}
The category of quasi-coherent modules on the Emerton-Gee moduli curve $\EG$ will be called the \emph{category of mod $p$ Langlands parameters with determinant $\omega\zeta$}, and denoted by $\LP_{\mathbf{\whG},0,\omega\zeta}$:
$$
\LP_{\mathbf{\whG},0,\omega\zeta}:=\QCoh(\EG).
$$
\end{Def*}

\subsection{A semisimple mod $p$ Langlands correspondence in families for $F=\bbQ_p$}

\begin{Pt*}  Let us consider $W$ to be a subgroup of $G$, by sending $s$ to the matrix 
$ \left(\begin{array}{cc}
0 & 1\\
1 & 0
\end{array} \right)
$
and by identifying the group $\Lambda$ with a subgroup of $T$ via $(1,0)\mapsto\diag(\varpi^{-1},1)$ and 
$(0,1)\mapsto\diag(1,\varpi^{-1})$. We obtain for example (recall that $u=(1,0)s\in W$)

$$
u =
\left (\begin{array}{cc}
0 & \varpi^{-1}\\
1 & 0
\end{array} \right),\hskip15pt
u^{-1} =
\left (\begin{array}{cc}
0 & 1\\
\varpi & 0
\end{array} \right),\hskip15pt
us =
\left (\begin{array}{cc}
\varpi^{-1} & 0\\
0 & 1
\end{array} \right),\hskip15pt
su =
\left (\begin{array}{cc}
1 & 0\\
0 & \varpi^{-1}
\end{array} \right).
$$
Moreover, $u^{2}=\diag(\varpi^{-1},\varpi^{-1})$.\footnote{Note that our element $u$ equals the element $u^{-1}$ in \cite{Be11},\cite{Br07} and \cite{V04}.} Since 
$$
\left (\begin{array}{cc}
0 & \varpi^{-1}\\
1 & 0
\end{array} \right) 
\left (\begin{array}{cc}
a & b\\
c & d
\end{array} \right)
\left (\begin{array}{cc}
0 & 1\\
\varpi & 0
\end{array} \right) = \left (\begin{array}{cc}
d & \varpi^{-1}c\\
\varpi b & a
\end{array} \right)
$$
the element $u\in G$ normalizes 
the group $I^{(1)}$.
\end{Pt*}

\begin{Pt*} \label{I(1)invariants}\label{some_invariants} Let $\Mod^{\rm sm}(k[G])$ be the category of smooth $G$-representations over $k$. 
Taking $I^{(1)}$-invariants yields a functor $\pi\mapsto \pi^{I^{(1)}}$ from $\Mod^{\rm sm} (k[G])$ to the category $\Mod(\cH^{(1)}_{\overline{\bbF}_q})$. If $F=\bbQ_p$, it induces a bijection between the irreducible $G$-representations and the irreducible
$\cH^{(1)}_{\overline{\bbF}_p}$-modules, under which supersingular representations correspond to supersingular Hecke modules \cite{V04}. 

\vskip5pt

For future reference, let us recall the $I^{(1)}$-invariants for some classes of representations. 
 If $\pi=\Ind_B^G(\chi)$ is a principal series representation with $\chi=\chi_1\otimes\chi_2$, then $\pi^{I^{(1)}}$ is a standard module in the component $\gamma:=\{\chi|_{\bbT}, \chi^s|_{\bbT}\}$. 
 
In the regular case, one chooses the ordering $(\chi|_{\bbT}, \chi^s|_{\bbT})$ on the set $\gamma$ and standard coordinates $x,y$. Then 
$$\Ind_B^G(\chi)^{I^{(1)}}=M(0,\chi(su),\chi (u^2),\chi|_{\bbT})=M(0,\chi_2(\varpi^{-1}),\chi_1(\varpi^{-1})\chi_2(\varpi^{-1}),\chi|_{\bbT})$$ 
In the non-regular case, one has 
$$\Ind_B^G(\chi)^{I^{(1)}}=M(\chi(su),\chi(u^2),\chi|_{\bbT})=M(\chi_2(\varpi^{-1}),\chi_1(\varpi^{-1})\chi_2(\varpi^{-1}),\chi|_{\bbT}).$$ 
These standard modules are irreducible if and only if $\chi\neq\chi^s$ \cite[4.2/4.3]{V04}.\footnote{Our formulas differ from \cite[4.2/4.3]{V04} by $\chi(\cdot)\leftrightarrow \chi(\cdot)^{-1}$, since we are working with left modules; also compare with the explicit
calculation with right convolution given in \cite[Appendix A.5]{V04}.}

\vskip5pt

Let $F=\bbQ_p$. If $\pi=\pi(r,0,\eta)$ is a standard supersingular representation with parameter $r=0,...,p-1$ and central character $\eta: \bbQ_p^{\times}\ra k^{\times}$, then $\pi^{I^{(1)}}$ is a supersingular module in the component $\gamma=\{\chi,\chi^s\}$ represented by the character $\chi:=(\omega^r\otimes 1)\cdot(\eta |_{\bbF_p^{\times}})$, cf. \cite[5.1/5.3]{Br07}. If $\pi$ is the trivial representation $\mathbbm{1}$ or the Steinberg representation {\rm St}, then $\gamma=1$ and $\pi^{I^{(1)}}$ is the character $(0,1)$ or $(-1,-1)$ respectively. 
\end{Pt*}

\begin{Pt*}\label{translate_the_action}
Let $\pi\in \Mod^{\rm sm} (k[G])$. Since $u\in G$ normalizes 
the group $I^{(1)}$, one has $I^{(1)}u I^{(1)}= u I^{(1)}$. It follows that the convolution action of the Hecke operator $U$ (resp. $U^2$) on $\pi^{I^{(1)}}$ is therefore induced by the action of $u$ (resp. $u^2$ on $\pi$). 
Similarly, the group $I^{(1)}$ is normalized by the Iwahori subgroup $I$ and $I/I^{(1)}\simeq \bbT$. It follows that
the convolution action of the operators $T_t, t\in \bbT$ on $\pi^{I^{(1)}}$ is the factorization of the $\mathbf{T}(o_F)$-action on $\pi$.

\end{Pt*}

\begin{Pt*}\label{ZGveeZvee}
We identify $F^{\times}$ with the center $Z(G)$ via $a\mapsto\diag(a,a)$. A (smooth) character 
$$
\zeta:  Z(G)=F^{\times}\lra k^\times
$$ 
is determined by its value $\zeta(\varpi^{-1})\in k^\times$ and its restriction $\zeta |_{o_F^{\times}}$. Since the latter is trivial on the subgroup $1+\varpi o_F$, we may view it as a character of $\bbF_q^{\times}$; we will write $\zeta |_{\bbF_q^{\times}}$ for this restriction in the following. Thus the group of characters of $Z(G)$ gets identified with the group of $k$-points of the group scheme $\cZ^{\vee} =(\bbF_q^{\times})^{\vee}\times \bbG_m$: 
$$Z(G)^\vee\iso \cZ^{\vee}(k),\; \zeta\mapsto (\zeta |_{\bbF_q^{\times}},\zeta(\varpi^{-1})).$$
\end{Pt*}

\begin{Lem*}\label{central_car_comp} 
Suppose that $\pi\in \Mod^{\rm sm} (k[G])$ has a central character $\zeta: Z(G)\rightarrow k^{\times}$. Then the Satake parameter $S(\pi^{I^{(1)}})$ of $\pi^{I^{(1)}}\in\Mod(\cH^{(1)}_{\overline{\bbF}_q})$ has central character $\zeta$, i.e. it is supported on the closed subscheme
$$
(V^{(1)}_{\mathbf{\whT},0}/W_0)_{(\zeta |_ {\bbF^{\times}_q} , \zeta(\varpi^{-1}))}\subset V^{(1)}_{\mathbf{\whT},0}/W_0.
$$
\end{Lem*}

\begin{proof} 
If $M$ is any $\cH^{(1)}_{\overline{\bbF}_q}$-module, then
$$
M=\bigoplus_{\gamma\in\bbT^{\vee}/W_0}\varepsilon_{\gamma}M=\bigoplus_{\gamma\in\bbT^{\vee}/W_0}\oplus_{\lambda\in\gamma}\varepsilon_{\lambda}M,
$$
and $\bbT\subset \overline{\bbF}_q[\bbT]\subset \cH^{(1)}_{\overline{\bbF}_q}$ acts on $\varepsilon_{\lambda}M$ through the character
$\lambda:\bbT\ra \bbF_q^{\times}$. Now if $M=\pi^{I^{(1)}}$, then the $\bbT$-action on $M$ is the factorization of the  $\mathbf{T}(o_F)$-action on $\pi$, cf. \ref{translate_the_action}. In particular, the restriction of the $\bbT$-action along the diagonal inclusion $\bbF_q^{\times}\subset\bbT$ is the factorization of the action of the central subgroup $o^{\times}_F\subset Z(G)$ on $\pi$, which is given by $\zeta|_{o^{\times}_F}$ by assumption. Hence
$$
\varepsilon_{\gamma}M\neq 0\quad\Longrightarrow\quad \forall\lambda\in\gamma,\ \lambda|_{\bbF_q^{\times}}=\zeta|_{\bbF_q^{\times}}\ \textrm{i.e.}\ \gamma|_{\bbF_q^{\times}}=\zeta|_{\bbF_q^{\times}}.
$$
Moreover, the element $u^2=\diag(\varpi^{-1},\varpi^{-1})\in Z(G)$ acts on $\pi$ by multiplication by $\zeta (\varpi^{-1})$ by assumption. Therefore, by \ref{translate_the_action}, the Hecke operator $z_2:=U^2\in \cH^{(1)}_{\overline{\bbF}_q}$ acts on $\pi^{I^{(1)}}$ by multiplication by $\zeta (\varpi^{-1})$. Thus we have obtained that $S(\pi^{I^{(1)}})$ is supported on 
$$
\coprod_{\gamma\in(\bbT^{\vee}/W_0)_{\reg},  \gamma |_ {\bbF^{\times}_q}= \zeta |_ {\bbF^{\times}_q} } V_{\mathbf{\whT},0,\zeta (\varpi^{-1})}\coprod_{\gamma\in(\bbT^{\vee}/W_0)_{\nonreg}, \gamma |_ {\bbF^{\times}_q}= \zeta |_ {\bbF^{\times}_q} } V_{\mathbf{\whT},0,\zeta (\varpi^{-1})}/W_0
=(V^{(1)}_{\mathbf{\whT},0}/W_0)_{(\zeta |_ {\bbF^{\times}_q} , \zeta(\varpi^{-1}))}.
$$
\end{proof}

\noindent Next, recall the twisting action of the group $\cZ^{\vee}(k)$ on the standard $\cH^{(1)}_{\overline{\bbF}_q}$-modules and their simple constituents \ref{twistHecke}. 

\begin{Prop*} \label{comp_twist} Let $\pi\in \Mod^{\rm ladm} (k[G])$ be irreducible or a reducible principal series representation. Let $\eta: F^{\times}\ra k^{\times}$ be a character.
Then $$(\pi\otimes\eta)^{I^{(1)}} =\pi^{I^{(1)}}. (\eta |_{\bbF_q^{\times}},\eta(\varpi^{-1}))$$
as $\cH^{(1)}_{\overline{\bbF}_q}$-modules.
 \end{Prop*}
 
 \begin{proof} An irreducible locally admissible representation, 
 being a finitely generated $k[G]$-module, is admissible \cite[2.2.19]{Em10}. 
 A principal series representation (irreducible or not) is always admissible \cite[4.1.7]{Em10}.
 The list of irreducible admissible smooth $G$-representations is given in \cite[Thm. 1.1]{H11b}. There are four families: principal series representations, supersingular representations, characters and twists of the Steinberg representation. 
 
 We first suppose that $\pi$ is a principal series representation (irreducible or not), i.e. of the form 
 $\Ind_B^G(\chi)$ with a character $\chi=\chi_1\otimes\chi_2$. Then $\pi\otimes\eta\simeq\Ind_B^G(\chi_1\eta\otimes\chi_2\eta)$.
 We use the results from \ref{I(1)invariants}. The modules $\pi^{I^{(1)}}$ and
 $(\pi\otimes\eta)^{I^{(1)}}$ are standard modules in the components
 $\gamma:=\{\chi|_{\bbT}, \chi^s|_{\bbT}\}$ and $\gamma (\eta |_{\bbF_q^{\times}})$ respectively. 
 Suppose that $\gamma$ is regular. We choose the ordering $(\chi|_{\bbT}, \chi^s|_{\bbT})$ and standard coordinates $x,y$. Then 
 $$\Ind_B^G(\chi)^{I^{(1)}}=M(0,\chi_2(\varpi^{-1}),\chi_1(\varpi^{-1})\chi_2(\varpi^{-1}),\chi|_{\bbT})$$
 and 
 $$\Ind_B^G(\chi_1\eta\otimes\chi_2\eta)^{I^{(1)}}=M(0,\chi_2(\varpi^{-1})\eta(\varpi^{-1}),\chi_1(\varpi^{-1})\chi_2(\varpi^{-1})\eta(\varpi^{-2}),(\chi|_{\bbT}).(\eta |_{\bbF_q^{\times}})).$$
 This shows $(\pi\otimes\eta)^{I^{(1)}} =\pi^{I^{(1)}}. (\eta |_{\bbF_q^{\times}},\eta(\varpi^{-1}))$ in the regular case.
 Suppose that $\gamma$ is non-regular.
  Then 
$$\Ind_B^G(\chi)^{I^{(1)}}=M(\chi_2(\varpi^{-1}),\chi_1(\varpi^{-1})\chi_2(\varpi^{-1}),\chi|_{\bbT})$$ 
and 
$$\Ind_B^G(\chi_1\eta\otimes\chi_2\eta)^{I^{(1)}}=
M(\chi_2(\varpi^{-1})\eta(\varpi^{-1}),\chi_1(\varpi^{-1})\chi_2(\varpi^{-1})\eta(\varpi^{-2}),(\chi|_{\bbT}).(\eta |_{\bbF_q^{\times}})).$$
This shows $(\pi\otimes\eta)^{I^{(1)}} =\pi^{I^{(1)}}. (\eta |_{\bbF_q^{\times}},\eta(\varpi^{-1}))$ in the non-regular case.

We now treat the case where $\pi$ is a character or a twist of the Steinberg representation.
Consider the exact sequence  
$$1\rightarrow\mathbbm{1} \rightarrow \Ind_B^G(1)\rightarrow {\rm St}\rightarrow 1.$$
According to \cite[4.4]{V04} the sequence of invariants 
$$ (S): 1\rightarrow\mathbbm{1}^{I^{(1)}} \rightarrow \Ind_B^G(1)^{I^{(1)}} \rightarrow {\rm St}^{I^{(1)}} \rightarrow 1$$ is still exact and $\mathbbm{1}^{I^{(1)}}$ resp. ${\rm St}^{I^{(1)}}$ is the trivial character $(0,1)$ resp. sign character $(-1,-1)$ in the Iwahori component $\gamma=1$.
Tensoring the first exact sequence with $\eta$ produces the exact sequence  
$$1\rightarrow\eta \rightarrow \Ind_B^G(1) \otimes \eta \rightarrow {\rm St}\otimes\eta\rightarrow 1.$$
Since the restriction $\eta |_{o_F^{\times} }$ is trivial on $1+\varpi o_F$, one has $(\eta\circ\det) |_ { I^{(1)}} = 1$ and so, as a sequence of $k$-vector spaces with $k$-linear maps, the sequence of invariants
 $$1\rightarrow\eta^{I^{(1)}} \rightarrow (\Ind_B^G(1)\otimes \eta )^{I^{(1)}} \rightarrow ({\rm St}\otimes\eta)^{I^{(1)}} \rightarrow 1$$ 
coincides with the sequence $(S)$. It is therefore an exact sequence of $\cH^{(1)}_{\overline{\bbF}_q}$-modules, with outer
terms being characters of $\cH^{(1)}_{\overline{\bbF}_q}$. From the discussion above, we deduce 
$$ (\Ind_B^G(1)\otimes \eta )^{I^{(1)}} =  \Ind_B^G(1)^{I^{(1)}}.( \eta |_{\bbF_q^{\times}},\eta(\varpi)^{-1})=
M(\eta(\varpi^{-1}), \eta(\varpi^{-2}),  1.(\eta |_{\bbF_q^{\times}})).$$
It follows then from \cite[1.1]{V04} that $\eta^{I^{(1)}}$ must be the trivial character $(0,\eta(\varpi^{-1}))$ in the component
 $1.(\eta |_{\bbF_q^{\times}})$ and  $({\rm St}\otimes\eta)^{I^{(1)}}$ must be the sign character $(-1,-\eta(\varpi^{-1}))$ in the component  $1.(\eta |_{\bbF_q^{\times}})$. This implies 
$$ \eta^{I^{(1)}}=\mathbbm{1}^{I^{(1)}}. ( \eta |_{\bbF_q^{\times}},\eta(\varpi)^{-1}) \hskip10pt \text{and} \hskip10pt
({\rm St}\otimes\eta)^{I^{(1)}} = {\rm St}^{I^{(1)}}. ( \eta |_{\bbF_q^{\times}},\eta(\varpi)^{-1}).$$
This proves the claim in the cases $\pi=\mathbbm{1}$ or $\pi={\rm St}$. 
If, more generally, $\pi=\eta'$ is a general character of $G$, then 
$$(\pi\otimes\eta)^{I^{(1)}}= (\eta'\eta)^{I^{(1)}}= \mathbbm{1}^{I^{(1)}}. ( (\eta'\eta) |_{\bbF_q^{\times}},(\eta'\eta)(\varpi)^{-1}) = \pi^{I^{(1)}}.(\eta |_{\bbF_q^{\times}},\eta(\varpi)^{-1} ).$$
On the other hand, if $\pi={\rm St}\otimes\eta'$ is a twist of Steinberg, then

$$(\pi\otimes\eta)^{I^{(1)}}= ({\rm St}\otimes (\eta'\eta))^{I^{(1)}}= {\rm St}^{I^{(1)}}. ( (\eta'\eta) |_{\bbF_q^{\times}},(\eta'\eta)(\varpi)^{-1}) = \pi^{I^{(1)}}.(\eta |_{\bbF_q^{\times}},\eta(\varpi)^{-1} ).$$

It remains to treat the case where $\pi$ is a supersingular representation. In this case $\pi\otimes\eta$ is also supersingular and the two modules $\pi^{I^{(1)}}$ and $(\pi\otimes\eta)^{I^{(1)}}$ are supersingular $\cH^{(1)}_{\overline{\bbF}_q}$-modules \cite[4.9]{V04}. Let $\gamma$ be the component of the module $\pi^{I^{(1)}}$. By \ref{translate_the_action}, the component of $(\pi\otimes\eta)^{I^{(1)}}$ equals $\gamma (\eta |_{\bbF_q^{\times}})$.
Moreover, if $U^2$ acts on $\pi^{I^{(1)}}$ via the scalar $z_2\in k^{\times}$, then 
$U^2$ acts on $(\pi\otimes\eta)^{I^{(1)}}$ via $z_2 (\eta\circ\det)(u^2)=z_2\eta(\varpi)^{-2}$, cf. \ref{translate_the_action}.
Since the supersingular modules are uniquely characterized by their component and their $U^2$-action,
we obtain $(\pi\otimes\eta)^{I^{(1)}}=\pi^{I^{(1)}}.(\eta |_{\bbF_q^{\times}},\eta(\varpi)^{-1} )$, as claimed. 
 \end{proof}
\vskip5pt

\begin{Pt*} \label{blocks}
Let $F=\bbQ_p$ with $p\geq 5$. We let  $\Mod_{\zeta}^{\rm l adm}(k[G])$ be the full subcategory of $\Mod^{\rm sm}(k[G])$ consisting of locally admissible representations having central character $\zeta$. By work of Pa\v{s}k\={u}nas \cite{Pas13}, the blocks $b$ of the category $\Mod_{\zeta}^{\rm l adm}(k[G])$, defined as certain equivalence classes of simple objects, can be parametrized by the set of isomorphism classes $[\rho]$ of semisimple continuous Galois representations $\rho : {\rm Gal}(\overline{\bbQ}_p/ \bbQ_p)\rightarrow \mathbf{\whG}(k)$ having determinant $\det \rho = \omega \zeta$, i.e. by the $k$-points of $\EG$. There are three types of blocks. Blocks of type $1$ are 
supersingular blocks. Each such block contains only one irreducible $G$-representation, which is supersingular. Blocks of type $2$ 
contain only two irreducible representations. These two representations are two generic principal series representations of the form
$\Ind_B^G(\chi_1\otimes\chi_2\omega^{-1})$ and $\Ind_B^G(\chi_2\otimes\chi_1\omega^{-1})$ (where $\chi_1\chi_2\neq 1,\omega^{\pm 1}$). 
There are four blocks of type $3$ which correspond to the four exceptional points. In the even case, each such block contains only three irreducible representations. These representations are of the form $\eta, {\rm St}\otimes\eta$ and $\Ind_B^G(\omega\otimes\omega^{-1}) \otimes\eta$. In the odd case, each block of type $3$ contains only one irreducible representation. It is of the form $\Ind_B^G(\chi\otimes\chi\omega^{-1})$. 
\end{Pt*}

\begin{Pt*} \label{ssLLC}
Let $F=\bbQ_p$ with $p\geq 5$.
Pa\v{s}k\={u}nas' parametrization $[\rho]\mapsto b_{[\rho]}$ is compatible with Breuil's semisimple mod $p$ local Langlands correspondence 
$$\rho\mapsto \pi(\rho)$$ for the group $G$ \cite{Br07,Be11}, in the sense that if $\rho$ has determinant $\omega\zeta$, then the simple constituents of the $G$-representation $\pi(\rho)$ lie in the block $b_{[\rho]}$ of $\Mod_{\zeta}^{\rm l adm}(k[G])$. 

The correspondence and 
the parametrizations (for varying $\zeta$) commute with twists: for a character $\eta: \bbQ_p^{\times}\ra k^{\times}$, $\pi(\rho\otimes\eta)=\pi(\rho)\otimes\eta$ and $b_{[\rho]}\otimes\eta=b_{[\rho\otimes\eta]}$. 
\end{Pt*}

\begin{Th*}\label{LLfamily}
Suppose $F=\bbQ_p$ with $p\geq 5$. Fix a character $\zeta:  Z(G)=\bbQ_p^{\times} \ra k^\times$, corresponding to a point $(\zeta|_{\bbF_p^{\times}},\zeta(p^{-1}))\in\cZ^{\vee}(k)$ under the identification $\cZ(G)^{\vee}\cong\cZ^{\vee}(k)$ from \ref{ZGveeZvee}. Let $(V^{(1)}_{\mathbf{\whT},0}/W_0)_{\zeta}$ be the space of mod $p$ Satake parameters with central character $\zeta$ and $X_{\zeta}$ be the moduli space of mod $p$ Langlands parameters with determinant $\omega\zeta$.
 
There exists a morphism of $k$-schemes  
$$ L_\zeta: (V^{(1)}_{\mathbf{\whT},0}/W_0)_{\zeta}\longrightarrow \EG$$ 
such that the quasi-coherent $\cO_{X_{\zeta}}$-module 
$$
L_{\zeta*}S(\cM_{\overline{\bbF}_p}^{(1)})|_{(V^{(1)}_{\mathbf{\whT},0}/W_0)_{\zeta}}
$$ 
equal to the push-forward along $L_{\zeta}$ of the restriction to $(V^{(1)}_{\mathbf{\whT},0}/W_0)_{\zeta}\subset V^{(1)}_{\mathbf{\whT},0}/W_0$ of the Satake parameter of the mod $p$ antispherical module $\cM_{\overline{\bbF}_p}^{(1)}$ interpolates the $I^{(1)}$-invariants of the semisimple mod $p$ Langlands correspondence
$$
\left.
\begin{array}{lllll}
X_{\zeta}(k) & \lra & \Mod_{\zeta}^{\rm l adm}(k[G]) & \lra & \Mod(\cH^{(1)}_{\overline{\bbF}_p}) \\
x & \lmapsto & \pi(\rho_x) & \lmapsto & \pi(\rho_x)^{I^{(1)}},
\end{array}
\right.
$$
in the sense that for all $x\in X_{\zeta}(k)$,
$$
\Big(\big(L_{\zeta*}S(\cM_{\overline{\bbF}_p}^{(1)})|_{(V^{(1)}_{\mathbf{\whT},0}/W_0)_{\zeta}}\big)\otimes_{\cO_{X_{\zeta}}}k(x)\Big)^{\sss}=\Big(\cM_{\overline{\bbF}_p}^{(1)}\otimes_{Z(\cH^{(1)}_{\overline{\bbF}_p})}(\sS_{\overline{\bbF}_p}^{(1)})^{-1}(\cO_{L_{\zeta}^{-1}(x)})\Big)^{\sss}\cong \pi(\rho_x)^{I^{(1)}}
$$
in $\Mod(\cH^{(1)}_{\overline{\bbF}_p})$.
\end{Th*}


\begin{Pt*}\label{prop_L_zeta}
The connected components of $(V^{(1)}_{\mathbf{\whT},0}/W_0)_{\zeta}$ are either regular and then of type $\bbA^1 \cup_0 \bbA^1$, or non-regular and then of type $\bbA^1$. The morphism $L_{\zeta}$ appearing in the theorem depends on the choice of an order of the two affine lines in each regular component. It is surjective and quasi-finite. Moreover, writing $L_{\zeta}^{\gamma}$ for its restriction to the connected component $(V^{\gamma}_{\mathbf{\whT},0}/W_0)_{\zeta}\subset (V^{(1)}_{\mathbf{\whT},0}/W_0)_{\zeta}$, one has:
\begin{itemize}
\item[(e)] \emph{Even case.} All connected components are of type $\bbA^1 \cup_0 \bbA^1$, except for the two `exterior' components which are of type $\bbA^1$. $L_{\zeta}^{\gamma}$ is an open immersion for any $\gamma$.

\item[(o)] \emph{Odd case.} All connected components are of type $\bbA^1 \cup_0 \bbA^1$.
$L_{\zeta}$ is an open immersion on all connected components, except for the two `exterior' ones. 
On an `exterior' component $\gamma$, the restriction of $L_{\zeta}^{\gamma}$ to one irreducible component $\bbA^1$ is an open immersion, and its restriction to the open complement $\bbG_m$ is a degree $2$ finite flat covering of its image, with branched locus equal to the intersection of this image with the exceptional locus of $X_{\zeta}$.
\end{itemize}
\end{Pt*}

\begin{Pt*} 
Note that the semisimple mod $p$ Langlands correspondence associates with any semisimple 
$\rho: {\rm Gal}(\overline{\bbQ}_p/ \bbQ_p)\rightarrow \mathbf{\whG}(k)$ a semisimple smooth $G$-representation $\pi(\rho)$ of length $1,2$ or $3$, hence whose semisimple $\cH^{(1)}_{\overline{\bbF}_p}$-module of $I^{(1)}$-invariants $\pi(\rho)^{I^{(1)}}$ has length $1,2$ or $3$. 
On the other hand, the antispherical map 
 $$
\xymatrix{
\ASph:  (V^{(1)}_{\mathbf{\whT},0}/W_0)(k)\ar[r] & \{\textrm{left $\cH^{(1)}_{\overline{\bbF}_q}$-modules}\}
}
$$
has an image consisting of $\cH^{(1)}_{\overline{\bbF}_q}$-modules are of length $1$ or $2$, cf. \cite[7.4.9]{PS} and  \cite[7.4.15]{PS}. Theorem \ref{LLfamily} combined with the properties \ref{prop_L_zeta} of the morphism $L_{\zeta}$ provide the following case-by-case elucidation of the $\cH^{(1)}_{\overline{\bbF}_p}$-modules  $\pi(\rho)^{ I^{(1)}}$.
\end{Pt*}

\begin{Cor*}\label{Morphism_V_to_X} 
Let $x\in \EG(k)$, corresponding to $\rho_x:{\rm Gal}(\overline{\bbQ}_p/ \bbQ_p)\rightarrow \mathbf{\whG}(k)$. Then the $\cH^{(1)}_{\overline{\bbF}_p}$-module  $\pi(\rho)^{ I^{(1)}}$ admits the following explicit description.
\begin{itemize}
\item[(i)]  If $x\in \EG^{irred}(k)$, then the fibre $L_\zeta^{-1}(x)=\{v\}$ has cardinality $1$ and 
$$
\pi(\rho_x)^{I^{(1)}}\simeq \ASph(v).
$$
It is irreducible and supersingular. 

\item[(ii)]  If $x\in \EG^{red}(k)\setminus \{ \text{the four exceptional points} \}$, then $L_\zeta^{-1}(x)=\{v_1,v_2\}$ has cardinality $2$ and 
$$
\pi(\rho_x)^{I^{(1)}}\simeq \ASph(v_1)\oplus \ASph(v_2).
$$
It has length $2$.

\item[(iiie)]  If $x\in \EG^{red}(k)$ is exceptional in the even case, then $L_\zeta^{-1}(x)=\{v_1,v_2\}$ has cardinality $2$ and 
$$
\pi(\rho_x)^{I^{(1)}}\simeq \ASph(v_1)^{\sss}\oplus \ASph(v_2).
$$
It has length $3$.

\item[(iiio)]  If $x\in \EG^{red}(k)$ is exceptional in the odd case, then $L_\zeta^{-1}(x)=\{v\}$ has cardinality $1$ and 
$$
\pi(\rho_x)^{I^{(1)}}\simeq \ASph(v)\oplus \ASph(v).
$$
It has length $2$.
\end{itemize}
\end{Cor*}

\begin{Pt*} \label{Lzetasets}
Now we proceed to the proof of \ref{LLfamily}, \ref{prop_L_zeta} and \ref{Morphism_V_to_X}. 

We start by defining the morphism $L_{\zeta}$ at the level of $k$-points. Let $v\in (V^{(1)}_{\mathbf{\whT},0}/W_0)_{\zeta}(k)$ and let its connected component be indexed by $\gamma\in\bbT^{\vee} /W_0$. \vskip5pt

1. Suppose that $\gamma$ is regular. Then
$\ASph(v)=\ASph^{\gamma}(v)$ is a simple two-dimensional $\cH^{\gamma}_{\overline{\bbF}_p}$-module, cf. \cite[7.4.9]{PS}.
Let $\pi\in \Mod^{\rm sm} (k[G])$ be the simple module, unique up to isomorphism, such that $\pi^{I^{(1)}}\simeq \ASph^{\gamma}(v)$, cf. \ref{I(1)invariants}. Then $\pi\in \Mod_{\zeta}^{\rm l adm}(k[G])$ with
$$
\zeta=(\zeta |_{\bbF_p^{\times}}, \zeta(p^{-1}))=(\gamma |_{\bbF_p^{\times}}, z_2)
$$
by \ref{central_car_comp}. Let $b$ be the block of $\Mod_{\zeta}^{\rm l adm}( k[G])$ which contains $\pi$. We define $L_{\zeta}(v)$ to be the point of $\EG(k)$ which corresponds to $b$.  

\vskip5pt 

2. Suppose that $\gamma$ is non-regular. 

(a) If $v\in D(2)_{\gamma}(k)$, then $\ASph(v)=\ASph^{\gamma}(2)(v)$ is 
a simple two-dimensional $\cH^{\gamma}_{\overline{\bbF}_p}$-module, cf. \cite[7.4.15]{PS}. As in the regular case, there is a simple module $\pi$, unique up to isomorphism, such that $\pi^{I^{(1)}}\simeq \ASph^{\gamma}(2)(v)$. It has central character $\zeta=(\gamma |_{\bbF_p^{\times}}, z_2)$ and there is a block $b$ of $\Mod_{\zeta}^{\rm l adm}( k[G])$ which contains $\pi$. We define $L_{\zeta}(v)$ to be the point of $\EG(k)$ which corresponds to $b$. 

(b) If  $v\in D(1)_{\gamma}(k)$, then $\ASph(v)^{\sss}$ is the direct sum of the two characters forming the antispherical pair 
$\ASph^{\gamma}(1)(v)=\{(0,z_1),(-1,-z_1)\}$ where $z_2=z^2_1$, cf. \cite[7.4.15]{PS}. As in the regular case, there are two simple modules $\pi_1$ and $\pi_2$, unique up to isomorphism, such that $\pi_1^{I^{(1)}}\simeq (0,z_1)$ and  $\pi_2^{I^{(1)}}\simeq (-1,-z_1)$ and $\pi_1, \pi_2$ have central character $\zeta=(\gamma |_{\bbF_p^{\times}}, z_2)$.
Moreover, we claim that there is a unique block $b$ of $\Mod_{\zeta}^{\rm l adm}( k[G])$ which contains both $\pi_1$ and $\pi_2$. 
Indeed, if $\gamma=\{ 1\otimes 1\}$ and $z_1=1$, then $\pi_1=\mathbbm{1}$ and $\pi_2={\rm St}$, cf. \ref{some_invariants}.
Then by \ref{comp_twist} it follows more generally that if $\gamma=\{ \omega^r\otimes\omega^r\}$, then $\pi_1=\eta$ and $\pi_2={\rm St}\otimes\eta$ with $\eta=(\eta |_{\bbF_p^{\times}},\eta(p^{-1})):=(\omega^r,z_1)$. Consequently $\pi_1,\pi_2$ are contained in a unique block $b$ of type $3$, cf. \ref{blocks}. We define $L_{\zeta}(v)$ to be the point of $\EG(k)$ which corresponds to $b$. 

\vskip5pt

Thus we have a well-defined map of sets $L_{\zeta}: (V^{(1)}_{\mathbf{\whT},0}/W_0)_{\zeta}(k)\longrightarrow \EG(k).$

\vskip5pt 

We show property (i) of \ref{Morphism_V_to_X}. Let $x\in\EG^{\rm irred}(k)$ and suppose $L_{\zeta}(v)=x$. 
Then $b_x$ is a supersingular block, contains a unique irreducible representation $\pi$, which is supersingular, and
$\pi=\pi(\rho_x)$, cf. \ref{blocks}-\ref{ssLLC}.
By definition of $L_{\zeta}$, one has $\ASph(v)\simeq \pi^{I^{(1)}}$. 
Since the antispherical map $\ASph$ is $1:1$ over supersingular modules, cf.  \cite[7.4.9]{PS} and
 \cite[7.4.15]{PS}, such a preimage $v$ of $x$ exists and is uniquely determined by $x$. 
Summarizing, we have $L_\zeta^{-1}(x)=\{v\}$ and $\ASph(v)\simeq \pi(\rho_x)^{I^{(1)}}$. This is property (i). 

\vskip5pt 

As a next step, we take a second character $\eta:\bbQ_p^{\times}\ra k^{\times}$ and show that the diagram 
$$
\xymatrix{
(V^{(1)}_{\mathbf{\whT},0}/W_0)_{\zeta}(k)\ar[r]^<<<<<<<{L_{\zeta}} \ar[d]^{\simeq}_{. \eta}  &  \EG(k) 
\ar[d]_{\simeq}^{(\cdot)\otimes\eta} \\
(V^{(1)}_{\mathbf{\whT},0}/W_0)_{\zeta\eta^2}(k)\ar[r]^<<<<{L_{\zeta\eta^2}}& X_{\zeta\eta^2}(k)
}
$$
commutes. Here, the vertical arrows are the bijections coming from \ref{twistSatake} and \ref{twistEG}. 
To verify the commutativity, let $v\in (V^{(1)}_{\mathbf{\whT},0}/W_0)_{\zeta}(k)$ and let its connected component be indexed by 
$\gamma\in\bbT^{\vee} /W_0$. Suppose that $\gamma$ is regular or that $\gamma$ is non-regular with $v\in D(2)_{\gamma}(k)$.
Let $\pi$ be the simple $G$-module with $\pi^{I^{(1)}}\simeq \ASph(v)$ and let $b_{[\rho]}$ be the block corresponding to the point 
$L_{\zeta}(v)$. By the equivariance property \ref{Asphequiv}, one has $\ASph(v.\eta)\simeq  \ASph(v).\eta$. Taking $I^{(1)}$-invariants
is compatible with twist, cf. \ref{comp_twist}, and so $L_{\zeta\eta^2}(v.\eta)$ corresponds to the block which contains the representation 
$\pi\otimes\eta$, i.e. to $b_{[\rho]}\otimes\eta=b_{[\rho\otimes\eta]}$, cf. \ref{ssLLC}, and so $L_{\zeta\eta^2}(v.\eta)=[\rho\otimes\eta]=L_{\zeta}(v).\eta$. 

If $v\in D(1)_{\gamma}(k)$, let $\pi_1$ and $\pi_2$ be the simple modules
such that $(\pi_1\oplus\pi_2)^{I^{(1)}}\simeq \ASph^{\gamma}(v)^{\sss}$. 
As before, we conclude from $\ASph(v.\eta)^{\sss}\simeq  \ASph(v)^{\sss}\otimes\eta$ 
that $L_{\zeta\eta^2}(v.\eta)$ corresponds to the block which contains
$\pi_1\otimes\eta$ and $\pi_2\otimes\eta$ and that $L_{\zeta\eta^2}(v.\eta)=L_{\zeta}(v).\eta$. The commutativity of the diagram is proved. 

\vskip5pt

Thus, we are reduced to prove that the map $L_{\zeta}$ comes from a morphism of $k$-schemes satisfying \ref{LLfamily} and the remaining parts of \ref{Morphism_V_to_X} in the two basic cases of a character $\zeta$ such that $\zeta(p^{-1})=1$ and $\zeta |_{\bbF_p^{\times}} \in \{1, \omega^{-1}\}$. This is established in the next two subsections. 
\end{Pt*}

\subsection{The morphism $L_{\zeta}$ in the basic even case}
Let $\zeta: \bbQ_p^{\times}\ra k^{\times}$ be the trivial character. Here we show that the map of sets $L_{\zeta}: (V^{(1)}_{\mathbf{\whT},0}/W_0)_{\zeta}(k)\ra\EG(k)$ that we have defined in \ref{Lzetasets} satisfies properties (ii) and (iiie) of \ref{Morphism_V_to_X}, and we define a morphism of $k$-schemes $L_{\zeta}: (V^{(1)}_{\mathbf{\whT},0}/W_0)_{\zeta}\longrightarrow \EG$ which coincides with the previous map of sets at the level of $k$-points. By construction, it will have the properties \ref{prop_L_zeta}. This will complete the proof of \ref{Morphism_V_to_X}, \ref{prop_L_zeta} and \ref{LLfamily} in the case of an even character. 

\begin{Pt*} \label{prop(ii)_even} We verify the properties (ii) and (iiie). We work over an irreducible component $\bbP^1$ with label 
$"\Sym^r\otimes\det^{a} \; | \;  \Sym^{p-3-r}\otimes\det^{r+1+a}"$ where 
$0\leq r\leq p-3$ and $0\leq a\leq p-2$, cf. \ref{structureEGeven}. 
On this component, we choose an affine coordinate $x$
around the double point having $\Sym^r\otimes\det^{a}$ as one of its Serre weights. Away from this point, we have $x\neq 0$ and 
the corresponding Galois representation has the form 

$$
\rho_x=
\left (\begin{array}{cc}
\unr(x)\omega^{r+1} & 0\\
0 & \unr(x^{-1})
\end{array} \right)\otimes\eta
$$
with $\eta=\omega^{a}$. By \cite[1.3]{Be11} or \cite[4.11]{Br07}, we have 
$$\pi(\rho_x)= \pi(r,x,\eta)^{\sss} \oplus  \pi([p-3-r],x^{-1},\omega^{r+1}\eta)^{\sss}=:\pi_1\oplus\pi_2$$ 
where $[p-3-r]$ denotes the unique integer in $\{0,...,p-2\}$ which is congruent to $p-3-r$ modulo $p-1$.
Now suppose that $L_{\zeta}(v)=x$. We distinguish two cases. 

\vskip5pt

1. {\it The generic case $0<r<p-3$.} In this case, the point $x$ lies on one of the `interior' components of the chain $\EG$, which has no exceptional points.  The length of $\pi(\rho_x)$ is $2$. Indeed, 
$\pi_1=\pi(r,x,\eta)$ and $\pi_2=\pi(p-3-r,x^{-1},\omega^{r+1}\eta)$ are two irreducible principal series representations \cite[Thm. 4.4]{Br07}. The block $b_x$
is of type $2$ and contains only these two irreducible representations, cf. \ref{blocks}-\ref{ssLLC}. We may write 
$$\pi_1=\Ind_B^G(\chi) \otimes\eta$$
with $\chi= \unr(x)\otimes\omega^{r}\unr(x^{-1})$, according to \cite[Rem. 4.4(ii)]{Br07}.
By our assumptions on $r$, the character $\chi |_{\bbT}= 1\otimes\omega^{r}$ is regular (i.e. different 
from its $s$-conjugate). We conclude from \ref{comp_twist} and \ref{some_invariants} that $\pi_1^{I^{(1)}}$ is a simple $2$-dimensional standard module in the regular component represented by the character $(1\otimes\omega^{r}).(\eta|_{\bbF_p^{\times}})= (\eta|_{\bbF_p^{\times}})\otimes(\eta|_{\bbF_p^{\times}})\omega^r\in\bbT^{\vee}$.
Similarly, we may write $$\pi_2=\Ind_B^G(\chi) \otimes\omega^{r+1}\eta$$
where now $\chi= \unr(x^{-1})\otimes\omega^{p-3-r}\unr(x)$.
By our assumptions on $r$, the character $\chi |_{\bbT}= 1 \otimes\omega^{p-3-r}$ is regular and we conclude, as above,
that the $I^{(1)}$-invariants $\pi_2^{I^{(1)}}$ form a simple $2$-dimensional standard module in the regular component 
represented by the character $(\eta|_{\bbF_p^{\times}})\omega^{r+1}\otimes(\eta|_{\bbF_p^{\times}})\omega^{r+1}\omega^{p-3-r}\in\bbT^{\vee}$. Note that the component of $\pi_1^{I^{(1)}}$ is different from the component of $\pi_2^{I^{(1)}}$, by our assumptions on $r$. 

We conclude from $L_{\zeta}(v)=x$ that either $\ASph(v)=\pi_1^{I^{(1)}}$ or $\ASph(v)=\pi_2^{I^{(1)}}$. Since for $\gamma$ regular, the map 
$\ASph^{\gamma}$ is a bijection
onto all simple $\cH^{\gamma}_{\overline{\bbF}_p}$-modules, cf.  \cite[7.4.9]{PS}, one finds that
$L_{\zeta}^{-1}(x)=\{v_1,v_2\}$ has cardinality $2$ and 
$$
\ASph(v_1)\oplus \ASph(v_2)\simeq \pi(\rho_x)^{I^{(1)}}.
$$
This settles property (ii) of \ref{Morphism_V_to_X} in the generic case. 

\vskip5pt

2. {\it The boundary cases $r\in \{0,p-3\}$.} In this case, the point $x$ lies on one of the two `exterior' components of $\EG$. On such a component, we will denote the variable $x$ rather by $z_1$, which is the notation\footnote{The reason for this notation will become clear in the discussion of the non-regular case in \ref{prop_mor_schemes_even}.} which we used already in \ref{structureEGeven}.

(a) Suppose that $z_1\neq \pm 1$. The length of $\pi(\rho_{z_1})$ is $2$. Indeed, as in the generic case, $\pi_1=\pi(r,z_1,\eta)$ and $\pi_2=\pi(p-3-r,z_1^{-1},\omega^{r+1}\eta)$ are two irreducible principal series representations. The block $b_{z_1}$ is of type $2$ and contains only these two irreducible representations. It follows, as above, that their invariants $\pi_1^{I^{(1)}}$ and $\pi_2^{I^{(1)}}$ are simple $2$-dimensional standard modules, in the components represented by $(\eta|_{\bbF_p^{\times}})\otimes(\eta|_{\bbF_p^{\times}})\omega^r\in\bbT^{\vee}$ and $(\eta|_{\bbF_p^{\times}})\omega^{r+1}\otimes(\eta|_{\bbF_p^{\times}})\omega^{r+1}\omega^{p-3-r}\in\bbT^{\vee}$ respectively. Since $r\in  \{0,p-3\}$, one of these components is regular, the other non-regular. In particular, the two components are different. We conclude from $L_{\zeta}(v)=z_1$ that either $\ASph(v)=\pi_1^{I^{(1)}}$ or $\ASph(v)=\pi_2^{I^{(1)}}$. Since for non-regular $\gamma$, the map $\ASph^{\gamma}(2)$ is a bijection from $D(2)_{\gamma}(k)$ onto all simple standard $\cH^{\gamma}_{\overline{\bbF}_p}$-modules, cf. \cite[7.4.15]{PS}, we may conclude as in the generic case: $L_{\zeta}^{-1}(z_1)=\{v_1,v_2\}$ has cardinality $2$ and 
$$
\ASph(v_1)\oplus \ASph(v_2)\simeq \pi(\rho_{z_1})^{I^{(1)}}.
$$
This settles property \ref{Morphism_V_to_X} (ii) in the remaining case $z_1\neq \pm 1$. 

(b) Suppose now that $z_1 = \pm 1$, i.e. we are at one of the four exceptional points. We will verify property (iiie).
The length of $\pi(\rho_{z_1})$ is $3$. Indeed, the representation $\pi(0,\pm 1,\eta)$ is a twist of the representation $\pi(0,1,1)$ (note that 
$\pi(r,z_1,\eta)\simeq \pi(r,-z_1,\unr(-1)\eta)$ according to \cite[Rem. 4.4(v)]{Br07}), 
which itself is an extension of $\mathbbm{1}$ by ${\rm St}$, cf. \cite[Thm. 4.4(iii)]{Br07}. As in the case (a), the representation $\pi_2=\pi(p-3,\pm 1,\omega\eta)$ is an irreducible principal series representation. The block $b_{z_1}$ is of type $3$ and contains only these three irreducible representations. The invariants $\pi_1^{I^{(1)}}$ form a direct sum of two antispherical characters in a non-regular component $\gamma$, whereas the invariants $\pi_2^{I^{(1)}}$ form a simple standard module in a regular component, as before.
Since for non-regular $\gamma$, the map $\ASph^{\gamma}(1)$ is a bijection from $D(1)_{\gamma}(k)$ onto 
all antispherical pairs of characters of $\cH^{\gamma}_{\overline{\bbF}_p}$, cf.  \cite[7.4.15]{PS}, 
we may conclude that $L_{\zeta}^{-1}(z_1)=\{v_1,v_2\}$ has cardinality $2$ with $v_1\in D(1)_{\gamma}(k)$ and $\ASph^{\gamma}(1)(v_1)^{\sss}=\pi_1^{I^{(1)}}$.
In particular,
$$
\ASph(v_1)^{\sss}\oplus \ASph(v_2)\simeq \pi(\rho_x)^{I^{(1)}}.
$$
This settles property \ref{Morphism_V_to_X} (iiie).
\end{Pt*}

\begin{Pt*} \label{prop_mor_schemes_even} We define a morphism of $k$-schemes $L_{\zeta}: (V^{(1)}_{\mathbf{\whT},0}/W_0)_{\zeta}\longrightarrow \EG$ which coincides on $k$-points with the map of sets $L_{\zeta}: (V^{(1)}_{\mathbf{\whT},0}/W_0)_{\zeta}(k)\longrightarrow \EG(k)$. We work over a connected component of $(V^{(1)}_{\mathbf{\whT},0}/W_0)_{\zeta}$, indexed by some $\gamma\in\bbT^{\vee}/W_0$. Let $v$ be a $k$-point of this component. 

Since $\zeta | _{ \bbF_p^{\times}} =1$, the connected components of $(V^{(1)}_{\mathbf{\whT},0}/W_0)_{\zeta}$ are indexed by the fibre $(\cdot)|_{\bbF_p^{\times}}^{-1}(1)$. This fibre consists of the $\frac{p-3}{2}$ regular components, represented by the characters of $\bbT$
$$\chi_k=\omega^{k}\otimes\omega^{-k}$$ for $k=1,...,\frac{p-3}{2}$, and of the two non-regular components, given by $\chi_0$ and $\chi_{\frac{p-1}{2}}$, cf. \ref{orbits}. 
We distinguish two cases. Note that $z_2=\zeta(p^{-1})=1$. 
\vskip5pt

1. {\it The regular case $0<k<\frac{p-1}{2}$.} We fix the order $\gamma=(\chi_k,\chi_k^s)$ on the set $\gamma$ and choose 
the standard coordinates $x,y$. According to  \cite[7.4.8]{PS}, our regular connected component identifies with two affine lines intersecting at the origin:
$$
V_{\mathbf{\whT},0,1}\simeq \bbA^1\cup_0 \bbA^1.
$$ 
Suppose that $v=(0,0)$ is the origin, so that $\ASph(v)$ is a supersingular module. Let $\pi(r,0,\eta)$ be the corresponding supersingular representation.
It corresponds to the irreducible Galois representation $ \rho(r,\eta)=\ind(\omega_2^{r+1})\otimes\eta$, in the notation of \cite[1.3]{Be11}, whence $L_{\zeta}(v)=[ \rho(r,\eta) ]$.
According to \ref{some_invariants}, the component of the Hecke module $\pi(r,0,\eta)^{I^{(1)}}$ is given by $(\omega^r\otimes 1)\cdot(\eta|_{\bbF_p^{\times}})$. Setting $\eta|_{\bbF_p^{\times}}  =\omega^{a}$,
this implies $(\omega^r\otimes 1)\cdot(\eta|_{\bbF_p^{\times}})=\omega^{r+a}\otimes\omega^{a}=\chi_k$ and hence $a=-k$ and $r=2k$. Therefore the Serre weights of the irreducible representation $\rho(r,\eta)$ are $\{ \Sym^{2k}\otimes\det ^{-k}, \Sym^{p-1-2k}\otimes\det ^{k}\}$, cf. \cite[1.9]{Br07}.

Comparing these pairs of Serre weights with the list 
\ref{structureEGeven} shows that the $\frac{p-3}{2}$ points $$\{ \text{origin $(0,0)$ on the component $(\chi_k,\chi_k^s)$} \}$$ 
for $0<k<\frac{p-1}{2}$ are mapped successively to the $\frac{p-3}{2}$ double points of the chain $\EG$. 

\vskip5pt Fix $0<k<\frac{p-1}{2}$ and consider the double point 
$$Q=L_{\zeta}( \text{origin $(0,0)$ on the component $\gamma=(\chi_k,\chi_k^s)$}).$$ As we have just seen, $Q$ lies on the irreducible 
component $\bbP^1$ whose label includes the weight $\Sym^{2k}\otimes\det^{-k}$ (i.e. on the component $" \Sym^{2k}\otimes\det^{-k} \; | \; \Sym^{p-3-2k}\otimes\det^{k+1}"$).
We fix an affine coordinate on this $\bbP^1$ around $Q$, which we will also call $x$ (there will be no risk of confusion with the standard coordinate above!).
 Away from $Q$, the affine coordinate $x\neq 0$ parametrizes Galois representations of the form 
$$
\rho_x=
\left (\begin{array}{cc}
\unr(x)\omega^{2k+1} & 0\\
0 & \unr(x^{-1})
\end{array} \right)\otimes\eta
$$
with $\eta:=\omega^{-k}$. As we have seen above,
$\pi(\rho_x)= \pi(2k,x,\eta) \oplus  \pi(p-3-2k,x^{-1},\omega^{r+1}\eta)=:\pi_1\oplus\pi_2$. 
Moreover, $\pi_1=\Ind_B^G(\chi) \otimes\eta$
with $\chi= \unr(x)\otimes\omega^{2k}\unr(x^{-1})$. Since 
$$(1\otimes\omega^{2k}).(\eta|_{\bbF_p^{\times}})=\omega^{-k}\otimes\omega^{k} = \chi_k^s\in\bbT^{\vee},$$ we deduce
 from the regular case of \ref{some_invariants} that
$$\pi_1^{I^{(1)}}=M(0,x,1,\chi_k^s)$$
 is a simple $2$-dimensional standard module. 
Note that $M(0,x,1,\chi_k^s)=M(x,0,1,\chi_k)$ according to \cite[Prop. 3.2]{V04}. 

Now suppose that $v=(x,0), x\neq 0,$ denotes
a point on the $x$-line of $\bbA_k ^1\cup_0 \bbA_k^1$. In particular, $\ASph^{\gamma}(v)= M(x,0,1,\chi_k)$. By our discussion, the point
$L_{\zeta}((x,0))$ corresponds to the block which contains $\pi_1$. Since $\pi_1$ lies in the block parametrized by $[\rho_x]$, cf. \ref{ssLLC}, it follows that
$$L_{\zeta}((x,0))=[\rho_x]=x\in \bbG_m \subset \bbP^1\subset X_{\zeta}.$$ Since $(0,0)$ maps to $Q$, i.e. to the point at $x=0$, the map
$L_{\zeta}$ identifies the whole affine $x$-line $\bbA^1=\{(x,0) : x\in k\}\subset V_{\mathbf{\whT},0,1}$
with the affine $x$-line $\bbA^1\subset \bbP^1\subset \EG$.

\vskip5pt

On the other hand, the double point $Q$ lies also on the irreducible 
component $\bbP^1$ whose labelling includes the other weight of $Q$, i.e. the weight $\Sym^{p-1-2k}\otimes\det^{k} $. 
We fix an affine coordinate $y$ on this $\bbP^1$ around $Q$. 
Away from $Q$, the coordinate $y\neq 0$ parametrizes Galois representations of the form 
$$
\rho_x=
\left (\begin{array}{cc}
\unr(y)\omega^{p-2k} & 0\\
0 & \unr(y^{-1})
\end{array} \right)\otimes\eta
$$
with $\eta:=\omega^{k}$. As in the first case, $\pi(\rho_y)$ contains 
$\pi_1:=  \pi(p-1-2k,y,\eta)= \Ind_B^G(\chi) \otimes\eta$
as a direct summand, where now $\chi= \unr(y)\otimes\omega^{p-1-2k}\unr(y^{-1})$. 
Since 
$$(1\otimes\omega^{p-1-2k}).(\eta|_{\bbF_p^{\times}})=\omega^k\otimes\omega^{-k}= \chi_k\in\bbT^{\vee},$$ we deduce, as above, that 
$\pi_1^{I^{(1)}}=M(0,y,1,\chi_k)$ is a simple $2$-dimensional standard module. 

Now suppose that $v=(0,y), y\neq 0,$ denotes
a point on the $y$-line of $\bbA_k ^1\cup_0 \bbA_k^1$. In particular, $\ASph^{\gamma}(v)= M(0,y,1,\chi_k)$.
By our discussion, the point $L_{\zeta}((0,y))$ corresponds to the block which contains $\pi_1$. Since $\pi_1$ lies in the block parametrized by 
$[\rho_y]$, cf. \ref{ssLLC},
it follows that 
$$L_{\zeta}((0,y))=[\rho_y]=y\in \bbG_m \subset \bbP^1\subset X_{\zeta}.$$ 
Since $(0,0)$ maps to $Q$, i.e. to the point at $y=0$, the map
$L_{\zeta}$ identifies the whole affine $y$-line  $\bbA^1=\{(0,y) : y\in k\} \subset V_{\mathbf{\whT},0,1}$
with the affine $y$-line $\bbA^1\subset \bbP^1\subset\EG$. 

\vskip5pt 

In this way, we get an open immersion of each regular connected component of $(V^{(1)}_{\mathbf{\whT},0}/W_0)_{\zeta}$ in the scheme $X_{\zeta}$, which coincides on $k$-points with the restriction of the map of sets $L_{\zeta}$.

\vskip5pt

2. {\it The non-regular case $k\in\{0,\frac{p-1}{2}\}$.} We choose the Steinberg coordinate $z_1$. 
According to  \cite[7.4.10]{PS}, our non-regular connected component identifies with an affine line :
$$
V_{\mathbf{\whT},0,z_2}/W_0\simeq \bbA^1.
$$ 
Suppose that $v=(0)$ is the origin, so that $\ASph(v)$ is a supersingular module. Let $\pi(r,0,\eta)$ be the corresponding supersingular representation
so that $L_{\zeta}(v)=[ \rho(r,\eta) ]$. Exactly as in the regular case, we may conclude that the Serre weights of the irreducible representation 
$\rho(r,\eta)$ are $\{ \Sym^{2k}\otimes\det ^{-k}, \Sym^{p-1-2k}\otimes\det ^{k}\}$. For the two values of $k=0$ and $k=\frac{p-1}{2}$ we find
 $\{ \Sym^{0}, \Sym^{p-1}\}$ and $\{ \Sym^{0}\otimes\det^{\frac{p-1}{2}}, \Sym^{p-1}\otimes\det ^{\frac{p-1}{2}}\}$ respectively. 
Comparing with the list 
\ref{structureEGeven} shows that the $2$ points $$\{ \text{origin $(0)$ on the component $(\chi_k=\chi_k^s)$} \}$$ for 
$k\in\{0,\frac{p-1}{2}\}$ are mapped to the $2$ smooth points in $\EG^{\rm irred}$, which lie on the two `exterior' components of $\EG$, cf.  \ref{structureEGeven}.

\vskip5pt Fix $k\in\{0,\frac{p-1}{2}\}$ and consider the point 
$$Q=L_{\zeta}( \text{origin $(0)$ on the component $\gamma=(\chi_k=\chi_k^s)$}).$$ As we have just seen, $Q$ lies on an `exterior' irreducible component $\bbP^1$ whose label includes the weight $\Sym^{0}\otimes\det^{k}$. 
We fix an affine coordinate on this $\bbP^1$ around $Q$, which we call $z_1$ (there will be no risk of confusion with the Steinberg coordinate above!).  Away from $Q$, the affine coordinate $z_1\neq 0$ parametrizes Galois representations of the form 
$$
\rho_{z_1}=
\left (\begin{array}{cc}
\unr(z_1)\omega & 0\\
0 & \unr(z_1^{-1})
\end{array} \right)\otimes\eta
$$
with $\eta:=\omega^{k}$. As in the regular case,
$\pi(\rho_{z_1})= \pi(0,z_1,\eta)^{\sss} \oplus  \pi(p-3,z_1^{-1},\omega\eta)^{\sss}$. 
Moreover, $\pi(0,z_1,\eta)=\Ind_B^G(\chi) \otimes\eta$
with $\chi= \unr(z_1)\otimes\unr(z_1^{-1})$ \footnote{The representations $\pi(0,z_1,\eta)$ constitute the {\it unramified} principal series of $G$.}.
Since 
$$(1\otimes 1).(\eta|_{\bbF_p^{\times}})=\omega^{k}\otimes\omega^{k} = \chi_k=\chi_k^s \in\bbT^{\vee},$$ we deduce
 from the non-regular case of \ref{some_invariants} that 
$ \pi(0,z_1,\eta)^{I^{(1)}}=M(z_1,1,\chi_k)$ is a $2$-dimensional standard module. 
Moreover, the standard module is simple if and only if $\chi\neq\chi^s$, i.e. if and only if 
$z_1\neq\pm 1$. 

Now let $v=z_1\neq 0$ denote
a nonzero point on our connected component $\bbA^1=V_{\mathbf{\whT},0,1}/W_0$.
Suppose that $z_1\neq\pm 1$, i.e. $v\in D(2)_{\gamma}$. In particular, $\ASph(v)= M(z_1,1,\gamma)$ is irreducible.
By our discussion, the point $L_{\zeta}(z_1)$ corresponds to the block (a block of type $2$) which contains $\pi(0,z_1,\eta)$. 
Suppose that $z_1=\pm 1$, i.e. $v\in D(1)_{\gamma}$. In particular,
$\ASph^{\sss}(v)=M(z_1,1,\chi_k)^{\sss}$ and again, $L_{\zeta}(z_1)$ corresponds to the block (now a block of type $3$) 
which contains the simple constituents of $\pi(0,z_1,\eta)^{\sss}$. In both cases, we conclude
$$L_{\zeta}(z_1)=[\rho_{z_1}]=z_1 \in \bbG_m \subset \bbP^1\subset\EG.$$
Since $(0)$ maps to $Q$, i.e. to the point at $z_1=0$, the map $L_{\zeta}$ identifies the whole 
$z_1$-line $\bbA^1=V_{\mathbf{\whT},0,1}/W_0 $
with the $z_1$-line $\bbA^1\subset \bbP^1\subset \EG$.

\vskip5pt

In this way, we get an open immersion of each non-regular connected component of $(V^{(1)}_{\mathbf{\whT},0}/W_0)_{\zeta}$ in the scheme $X_{\zeta}$, which coincides on $k$-points with the restriction of the map of sets $L_{\zeta}$.   
\end{Pt*}

\subsection{The morphism $L_{\zeta}$ in the basic odd case}
Let $\zeta:=\omega^{-1}: \bbQ_p^{\times}\ra k^{\times}$. Here we show that the map of sets $L_{\zeta}: (V^{(1)}_{\mathbf{\whT},0}/W_0)_{\zeta}(k)\ra\EG(k)$ that we have defined in \ref{Lzetasets} satisfies properties (ii) and (iiio) of \ref{Morphism_V_to_X}, and we define a morphism of $k$-schemes $L_{\zeta}: (V^{(1)}_{\mathbf{\whT},0}/W_0)_{\zeta}\longrightarrow \EG$ which coincides with the previous map of sets at the level of $k$-points. By construction, it will have the properties \ref{prop_L_zeta}.  This will complete the proof of \ref{Morphism_V_to_X}, \ref{prop_L_zeta} and \ref{LLfamily} in the case of an odd character.

\begin{Pt*}  \label{prop(ii)_odd} We verify properties (ii) and (iiio). We work over an irreducible component $\bbP^1$ with label \newline 
$"\Sym^r\otimes\det^{a} \; | \;  \Sym^{p-3-r}\otimes\det^{r+1+a}"$ where 
$1\leq r\leq p-2$ and $0\leq a\leq p-2$, cf. \ref{structureEGodd}. We distinguish two cases.

\vskip5pt

1. {\it The generic case $r \neq p-2$.} In this case, the irreducible component  of $\EG$ we consider is an `interior' component and has no exceptional points. On this component, we choose an affine coordinate $x$ around the double point having $\Sym^r\otimes\det^{a}$ as one of its Serre weights. Away from this point, we have $x\neq 0$ and the corresponding Galois representation has the form 

$$
\rho_x=
\left (\begin{array}{cc}
\unr(x)\omega^{r+1} & 0\\
0 & \unr(x^{-1})
\end{array} \right)\otimes\eta
$$
with $\eta=\omega^{a}$. As before, we have
$$\pi(\rho_x)= \pi(r,x,\eta)^{\sss} \oplus  \pi([p-3-r],x^{-1},\omega^{r+1}\eta)^{\sss}.$$ 
The length of $\pi(\rho_x)$ is $2$. Indeed, by our assumptions on $r$, the principal series representations
$\pi(r,x,\eta)$ and $\pi(p-3-r,x^{-1},\omega^{r+1}\eta)$ are irreducible and the block $b_x$ contains only these two irreducible representations. We may follow the argument of the generic case of \ref{prop(ii)_even} word for word and deduce property \ref{Morphism_V_to_X} (ii). 

\vskip5pt

2. {\it The two boundary cases $r=p-2$.} In this case, the irreducible component is one of the two `exterior' components with labels
$"\Sym^{p-2}  \;  | \;  "\Sym^{-1}" "$ or $""\Sym^{-1}\det^{\frac{p-1}{2}}" \;  | \;  \Sym^{p-2}\det^{\frac{p-1}{2}} "$. Points of 
the open locus $\EG^{\rm red}$ lying on such a component correspond to twists of unramified Galois representations of the form 
$$
\rho_{x+x^{-1}}=
\left (\begin{array}{cc}
\unr(x) & 0\\
0 & \unr(x^{-1})
\end{array} \right)\otimes\eta
$$
with $\eta=1$ or $\eta= \omega^{\frac{p-1}{2}}$. Let us concentrate on one of the two components, i.e. let us fix $\eta$. 

Mapping an unramified Galois representation $\rho_{x+x^{-1}}$ to $t:=x+x^{-1}\in k$ identifies this open locus with
the $t$-line $\bbA^1\subset \bbP^1$. We have $$\pi(\rho_{t})= \pi(p-2,x,\eta)^{\sss}\oplus  \pi(p-2,x^{-1},\eta)^{\sss}=:\pi_1\oplus\pi_2$$
since $[p-3-(p-2)]=p-2$ (indeed, $p-3-(p-2)=-1\equiv p-2\mod (p-1)$). The length of $\pi(\rho_t)$ is $2$. Indeed, 
$\pi_1=\pi(p-2,x,\eta)$ and $\pi_2=\pi(p-2,x^{-1},\eta)$ are two irreducible principal series representations 
and the block $b_{t}$ contains only these two irreducible representations. They are isomorphic if and only if $x=\pm 1$, i.e. if and only if $t=\pm 2$ is an exceptional point. In this case, $b_t$ contains only one irreducible representation and is of type $3$, otherwise it is of type 2.

We may write 
$$\pi_1=\Ind_B^G(\chi) \otimes\eta$$
with $\chi= \unr(x)\otimes\omega^{p-2}\unr(x^{-1})$. Similarly for $\pi_2$. 
The character $\chi |_{\bbF_p^{\times}}= 1\otimes\omega^{p-2}$ is regular (i.e. different 
from its $s$-conjugate) and we are in the regular case of \ref{some_invariants}. 
We conclude that $\pi_{1}^{I^{(1)}}= M(0,x,1,(1\otimes\omega^{p-2}).\eta)$ and $\pi_{2}^{I^{(1)}}= M(0,x^{-1},1,(1\otimes\omega^{p-2}).\eta)$ are both simple $2$-dimensional standard modules in the regular component $\gamma$ represented by the character 
$(1\otimes\omega^{p-2}).(\eta|_{\bbF_p^{\times}})= (\eta|_{\bbF_p^{\times}})\otimes(\eta|_{\bbF_p^{\times}})\omega^{p-2}\in\bbT^{\vee}$. They are isomorphic if and only if $t=\pm 2$. We choose an order $\gamma=(  (\eta|_{\bbF_p^{\times}})\otimes(\eta|_{\bbF_p^{\times}})\omega^{p-2},  (\eta|_{\bbF_p^{\times}})\omega^{p-2}\otimes(\eta|_{\bbF_p^{\times}}))$ on the set $\gamma$. Then from $L_{\zeta}(v)=t$ we get that either 
$\ASph^{\gamma}(v)=\pi_1^{I^{(1)}}$ or $\ASph^{\gamma}(v)=\pi_2^{I^{(1)}}$. 
Since for regular $\gamma$, the map $\ASph^{\gamma}$ is a bijection
onto all simple $\cH^{\gamma}_{\overline{\bbF}_p}$-modules, cf.  \cite[7.4.9]{PS}, one finds that
$L_{\zeta}^{-1}(t)=\{v_1,v_2\}$ has cardinality $2$ if $t\neq \pm 2$ and then  
$$
\ASph(v_1)\oplus \ASph(v_2)\simeq \pi(\rho_t)^{I^{(1)}}.
$$
This settles property \ref{Morphism_V_to_X} (ii). In turn, if $t=\pm 2$ is an exceptional point, then $L_{\zeta}^{-1}(t)=\{v\}$ has cardinality $1$ and 
 $$
\ASph(v)\oplus \ASph(v)\simeq \pi(\rho_t)^{I^{(1)}}.
$$
This settles property \ref{Morphism_V_to_X} (iiio).
\end{Pt*}

\begin{Pt*}  \label{prop_mor_schemes_odd} 
We define a morphism of $k$-schemes $L_{\zeta}: (V^{(1)}_{\mathbf{\whT},0}/W_0)_{\zeta}\longrightarrow \EG$ which coincides on $k$-points with the map of sets $L_{\zeta}: (V^{(1)}_{\mathbf{\whT},0}/W_0)_{\zeta}(k)\longrightarrow \EG(k)$. We work over a connected component of $(V^{(1)}_{\mathbf{\whT},0}/W_0)_{\zeta}$, indexed by some $\gamma\in\bbT^{\vee}/W_0$. Let $v$ be a $k$-point of this component. 

Since $\zeta | _{ \bbF_p^{\times}} =\omega^{-1}$, the connected components of $(V^{(1)}_{\mathbf{\whT},0}/W_0)_{\zeta}$ are indexed by the fibre $(\cdot)|_{\bbF_p^{\times}}^{-1}(\omega^{-1})$. This fibre consists of the $\frac{p-1}{2}$ regular components, represented by the characters 
$$\chi_k=\omega^{k-1}\otimes\omega^{-k}$$ for $k=1,...,\frac{p-1}{2}$, cf. \ref{orbits}. Recall that $z_2=\zeta(p)=1$. 
\vskip5pt

Fix an order $\gamma=(\chi_k,\chi_k^s)$ on the set $\gamma$ and choose 
standard coordinates $x,y$. According to  \cite[7.4.8]{PS}, our regular connected component identifies with two affine lines intersecting at the origin:
$$
V_{\mathbf{\whT},0,1}\simeq \bbA^1\cup_0 \bbA^1.
$$ 
Suppose that $v=(0,0)$ is the origin, so that $\ASph(v)$ is a supersingular module. Let $\pi(r,0,\eta)$ be the corresponding supersingular representation.
It corresponds to the irreducible Galois representation $ \rho(r,\eta)$, in the notation of \cite[1.3]{Be11}, whence $L_{\zeta}(v)=[ \rho(r,\eta) ]$.
According to \ref{some_invariants}, the component of $\pi(r,0,\eta)^{I^{(1)}}$ is given by $(\omega^r\otimes 1)\cdot (\eta |_{\bbF^{\times}_p})$. Setting $\eta |_{\bbF^{\times}_p}=\omega^{a}$, this implies
$(\omega^r\otimes 1)\cdot(\eta |_{\bbF^{\times}_p})=\omega^{r+a}\otimes\omega^{a}=\chi_k$ and hence $a=-k$ and $r=2k-1$. The Serre weights of the irreducible representation $\rho(r,\eta)$ are therefore $\{ \Sym^{2k-1}\otimes\det ^{-k}, \Sym^{p-2k}\otimes\det ^{k-1}\}$, cf. \cite[1.9]{Br07}. 

Comparing these pairs of Serre weights with the list 
\ref{structureEGodd} shows that the $\frac{p-1}{2}$ points $$\{ \text{origin $(0,0)$ on the component $(\chi_k,\chi_k^s)$} \}$$ 
for $k=1,..., \frac{p-1}{2}$ are mapped successively to the $\frac{p-1}{2}$ double points of the chain $\EG$. We distinguish two cases. 

\vskip5pt 

 {\it 1. The generic case $1<k<\frac{p-1}{2}$}.  In this case, the argument proceeds as in the regular case of \ref{prop_mor_schemes_even}. Consider the double point $$Q=L_{\zeta}( \text{origin $(0,0)$ on the component $\gamma=(\chi_k,\chi_k^s)$}).$$ As we have just seen, $Q$ lies on an `interior' irreducible component $\bbP^1$ whose label includes the weight $\Sym^{2k-1}\otimes\det^{-k}$. 
We fix an affine coordinate on this $\bbP^1$ around $Q$, which we will also call $x$. Away from $Q$, the affine coordinate $x\neq 0$ parametrizes Galois representations of the form 
$$
\rho_x=
\left (\begin{array}{cc}
\unr(x)\omega^{2k} & 0\\
0 & \unr(x^{-1})
\end{array} \right)\otimes\eta
$$
with $\eta:=\omega^{-k}$. As we have seen above,
$\pi(\rho_x)= \pi(2k-1,x,\eta) \oplus  \pi(p-3-2k+1,x^{-1},\omega^{2k}\eta)=:\pi_1\oplus\pi_2$. 
Moreover, $\pi_1=\Ind_B^G(\chi) \otimes\eta$
with $\chi= \unr(x)\otimes\omega^{2k-1}\unr(x^{-1})$. Since 
$$(1\otimes\omega^{2k-1}).(\eta|_{\bbF_p^{\times}})=\omega^{-k}\otimes\omega^{k-1} = \chi_k^s\in\bbT^{\vee},$$ we deduce
 from the regular case of \ref{some_invariants} that 
$\pi_1^{I^{(1)}}=M(0,x,1,\chi_k^s)$ is a simple $2$-dimensional standard module. 
Note that $M(0,x,1,\chi_k^s)=M(x,0,1,\chi_k)$ according to \cite[Prop. 3.2]{V04}. 

Now suppose that $v=(x,0)$, $x\neq 0$, denotes a nonzero
point on the $x$-line of $\bbA^1\cup_0 \bbA^1$. In particular, $\ASph^{\gamma}(v)= M(x,0,1,\chi_k)$. Our discussion
shows that the point  $L_{\zeta}((x,0))$ corresponds to the block which contains $\pi_1$. Since $\pi_1$ lies in the block parametrized by $[\rho_x]$, cf. \ref{ssLLC}, it follows that
$$L_{\zeta}((x,0))=[\rho_x]=x\in \bbG_m \subset \bbP^1.$$ Since $(0,0)$ maps to $Q$, i.e. to the point at $x=0$, the map
$L_{\zeta}$ identifies the whole affine $x$-line $\bbA^1=\{(x,0) : x\in k\}\subset V_{\mathbf{\whT},0,1}$
with the affine $x$-line $\bbA^1\subset \bbP^1\subset \EG$.

\vskip5pt

On the other hand, the double point $Q$ also lies on the irreducible 
component whose labelling includes the other weight of $Q$, i.e. the weight $\Sym^{p-2k}\otimes\det^{k-1} $. 
We fix an affine coordinate $y$ on this $\bbP^1$ around $Q$. 
Away from $Q$, the coordinate $y\neq 0$ parametrizes Galois representations of the form 
$$
\rho_y=
\left (\begin{array}{cc}
\unr(y)\omega^{p-2k+1} & 0\\
0 & \unr(y^{-1})
\end{array} \right)\otimes\eta
$$
with $\eta:=\omega^{k-1}$. As in the first case, $\pi(\rho_y)$ contains 
$\pi_1:=  \pi(p-2k,y,\eta)= \Ind_B^G(\chi) \otimes\eta$
as a direct summand, where now $\chi= \unr(y)\otimes\omega^{p-2k}\unr(y^{-1})$. 
Since 
$$(1\otimes\omega^{p-2k}).(\eta|_{\bbF_p^{\times}})=\omega^{k-1}\otimes\omega^{-k}= \chi_k\in\bbT^{\vee},$$ we deduce
 from the regular case of \ref{some_invariants} that 
$\pi_1^{I^{(1)}}=M(0,y,1,\chi_k)$ is a simple $2$-dimensional standard module. 

Now suppose that $v=(0,y)$, $y\neq 0,$ denotes
a nonzero point on the $y$-line of $\bbA^1\cup_0 \bbA^1$. In particular, $\ASph^{\gamma}(v)= M(0,y,1,\chi_k)$.
Our discussion shows that the point $L_{\zeta}((0,y))$ corresponds to the block which contains $\pi_1$, parametrized by  $[\rho_y]$.
Hence
$$L_{\zeta}((0,y))=[\rho_y]=y\in \bbG_m \subset \bbP^1.$$ Since $(0,0)$ maps to $Q$, i.e. to the point at $y=0$, the map
$L_{\zeta}$ identifies the whole $y$-line  $\bbA^1=\{(0,y) : y\in k\} \subset V_{\mathbf{\whT},0,1}$
with the affine $y$-line $\bbA^1\subset \bbP^1\subset\EG$. 

\vskip5pt 

In this way, we get an open immersion of each connected component $(V^{\gamma}_{\mathbf{\whT},0}/W_0)_{\zeta}$ of $(V^{(1)}_{\mathbf{\whT},0}/W_0)_{\zeta}$ such that $\gamma=(\chi_k,\chi_k^s)$ with $1<k<\frac{p-1}{2}$, in the scheme $X_{\zeta}$, which coincides on $k$-points with the restriction of the map of sets $L_{\zeta}$.

\vskip5pt

 {\it 2. The two boundary cases $k\in \{ 1,\frac{p-1}{2}\}$}. Consider the double point $$Q=L_{\zeta}( \text{origin $(0,0)$ on the component $\gamma=(\chi_k,\chi_k^s)$}).$$ As we have just seen, $Q$ lies on an `interior' irreducible component $\bbP^1$ 
 whose label includes the weight $\Sym^{1}\otimes\det^{-1}$ (for $k=1$) or the weight $\Sym^{1}\otimes\det^{\frac{p-3}{2}}$ (for $k=\frac{p-1}{2}$).
 We fix an affine coordinate on this $\bbP^1$ around $Q$, which we will call $z$.
 Away from $Q$, the coordinate $z\neq 0$ parametrizes Galois representations of the form 
$$
\rho_z=
\left (\begin{array}{cc}
\unr(z)\omega^{2} & 0\\
0 & \unr(z^{-1})
\end{array} \right)\otimes\eta
$$
with $\eta=\omega^{-1}$ or $\eta=\omega^{\frac{p-3}{2}}$. 

Let $k=1$, i.e. $\eta=\omega^{-1}$. 
Following the argument in the generic case word for word, we may conclude that 
$L_{\zeta}$ identifies the $x$-line $\bbA^1=\{(x,0) : x\in k\}\subset V_{\mathbf{\whT},0,1}$
with the $z$-line $\bbA^1\subset \bbP^1\subset \EG$. 

Let $k=\frac{p-1}{2}$, i.e. $\eta=\omega^{\frac{p-3}{2}}$. 
As in the generic case, we may conclude that 
$L_{\zeta}$ identifies the $y$-line $\bbA^1=\{(0,y) : y\in k\}\subset V_{\mathbf{\whT},0,1}$
with the $z$-line $\bbA^1\subset \bbP^1\subset \EG$. 

\vskip5pt

On the other hand, the double point $Q$ lies also on the irreducible 
component $\bbP^1$ whose labelling includes the other weight of $Q$, i.e. the weight $\Sym^{p-2} $ (for $k=1$) or the weight 
$\Sym^{p-2}\otimes\det^{\frac{p-1}{2}} $ (for $k=\frac{p-1}{2}$). These are the two `exterior' components.
 Points of the open locus $\EG^{\rm red}$ lying on such a component correspond to unramified (up to twist) Galois representations of the form 
$$
\rho_{t}=
\left (\begin{array}{cc}
\unr(z) & 0\\
0 & \unr(z^{-1})
\end{array} \right)\otimes\eta
$$
where $\eta=1$ (for $k=1$) or $\eta= \omega^{\frac{p-1}{2}}$  (for $k=\frac{p-1}{2}$) and with $t=z+z^{-1}\in\bbA^1\subset\bbP^1$. 
 As in the boundary case of \ref{prop(ii)_odd}, we have $\pi(\rho_t)= \pi(p-2,z,\eta)\oplus  \pi(p-2,z^{-1},\eta)=:\pi_1\oplus\pi_2$
and these are irreducible principal series representations. We may write $\pi_1=\Ind_B^G(\chi)\otimes\eta$
with $\chi= \unr(z)\otimes\omega^{p-2}\unr(z^{-1})$. The character $\chi |_{\bbF_p^{\times}}= 1\otimes\omega^{p-2}$ is regular (i.e. different 
from its $s$-conjugate) and we are in the regular case of \ref{some_invariants}. 
We conclude that $$\pi_{1}^{I^{(1)}}= M(0,z,1,(1\otimes\omega^{p-2}).\eta)$$ is a simple $2$-dimensional standard module in the regular component 
represented by the character $$(1\otimes\omega^{p-2}).(\eta|_{\bbF_p^{\times}})= (\eta|_{\bbF_p^{\times}})\otimes(\eta|_{\bbF_p^{\times}})\omega^{p-2}=(\eta|_{\bbF_p^{\times}})\otimes(\eta|_{\bbF_p^{\times}})\omega^{-1}\in\bbT^{\vee}.$$
This latter character equals $\chi_1$ for $\eta=1$ and $(\chi_{\frac{p-1}{2}})^s$ for $\eta=\omega^{\frac{p-1}{2}}$
(indeed, note that $\frac{p-1}{2}\equiv - \frac{p-1}{2}\mod p-1$). 

Now suppose that $k=1$, i.e. $\eta=1$. Let $v=(0,y)$, $y\neq 0$, be a nonzero point on the $y$-line of $\bbA^1\cup_0 \bbA^1$. In particular, 
$\ASph^{\gamma}(v)= M(0,y,1,\chi_1)$. Our discussion shows that the point $L_{\zeta}((0,y))$ corresponds to the block which contains 
$\pi_{1}$, i.e. which is parametrized by $[\rho_t]$. It follows that $$L_{\zeta}((0,y))=[\rho_t]=t=y+y^{-1}\in \bbA^1 \subset \bbP^1.$$ Since $(0,0)$ maps to $Q$, i.e. to the point at $t=\infty$, the map of sets
$L_{\zeta}$ maps the $k$-points of the whole affine $y$-line  $\bbA^1=\{(0,y) : y\in k\} \subset V_{\mathbf{\whT},0,1}$
to the $k$-points of the whole `left exterior' component $\bbP^1\subset\EG$ via the formula
\begin{eqnarray*}
\bbA^1 & \lra & \bbP^1\\
y &\lmapsto & 
\left\{ \begin{array}{ll} 
y+y^{-1} & \textrm{if $y\neq 0$} \\
\infty=Q & \textrm{if $y=0$}.
\end{array}
\right.
\end{eqnarray*}
This formula is algebraic: indeed, for $y\in\bbA^1\setminus \{ \pm i \}$ (where $\pm i$ are the roots of the polynomial
$f(y)=y^2+1$), we have $y+y^{-1}\neq 0$ and $(y+y^{-1})^{-1}=y/(y^2+1)$, which is equal to $0$ at $y=0$. Moreover, it glues at the origin $(0,0)$ with the open immersion of the $x$-line of $V_{\mathbf{\whT},0,1}=\bbA^1\cup_0 \bbA^1$ in $X_{\zeta}$ defined above, since both map $(0,0)$ to $Q$. We take the resulting morphism of $k$-schemes 
$\bbA^1\cup_0 \bbA^1\ra X_{\zeta}$ as the definition of $L_{\zeta}$ on the connected component $(V^{(\chi_1,\chi_1^s)}_{\mathbf{\whT},0}/W_0)_{\zeta}$ of $(V^{(1)}_{\mathbf{\whT},0}/W_0)_{\zeta}$. Note that its restriction to the open subset $\{y\neq 0\}$ in the $y$-line $\bbA^1$ is the morphism $\bbG_m\ra\bbA^1$ corresponding to the ring extension
$$
k[t] \lra k[y,y^{-1}]=k[t][y]/(y^2-ty+1),
$$
and that the discriminant $t^2-4$ of $y^2-ty+1\in k[t][y]$ vanishes precisely at the two exceptional points $t=\pm2$.

Suppose $k=\frac{p-1}{2}$, i.e. $\eta=\omega^{\frac{p-1}{2}}$. Let $v=(x,0)$, $x\neq 0$, denote
a nonzero point on the $x$-line of $\bbA^1\cup_0 \bbA^1$. In particular, 
$$\ASph^{\gamma}(v)= M(0,x,1,(\chi_{\frac{p-1}{2}})^s)=M(x,0,1,\chi_{\frac{p-1}{2}}).$$ 
Our discussion shows that the point $L_{\zeta}((x,0))$ corresponds to the block which contains $\pi_{1}$, i.e. which is parametrized by $[\rho_t]$. 
It follows that $L_{\zeta}((x,0))=[\rho_t]=t=x+x^{-1}\in \bbA^1 \subset \bbP^1$. Since $(0,0)$ maps to the point $Q$ at $t=\infty$, 
the map of sets $L_{\zeta}$ maps the $k$-points of the whole affine $x$-line  $\bbA^1=\{(x,0) : y\in k\} \subset V_{\mathbf{\whT},0,1}$
to the $k$-points of the whole `right exterior' component $\bbP^1\subset\EG$ via the formula
\begin{eqnarray*}
\bbA^1 & \lra & \bbP^1\\
x &\lmapsto & 
\left\{ \begin{array}{ll} 
x+x^{-1} & \textrm{if $x\neq 0$} \\
\infty=Q & \textrm{if $x=0$}.
\end{array}
\right.
\end{eqnarray*}
This formula is algebraic. Moreover, it glues at the origin $(0,0)$ with the open immersion of the $y$-line of $V_{\mathbf{\whT},0,1}=\bbA^1\cup_0 \bbA^1$ in $X_{\zeta}$ defined above, since both map $(0,0)$ to $Q$. We take the resulting morphism of $k$-schemes 
$\bbA^1\cup_0 \bbA^1\ra X_{\zeta}$ as the definition of $L_{\zeta}$ on the connected component $(V^{(\chi_{\frac{p-1}{2}},(\chi_{\frac{p-1}{2}})^s)}_{\mathbf{\whT},0}/W_0)_{\zeta}$ of $(V^{(1)}_{\mathbf{\whT},0}/W_0)_{\zeta}$.
\end{Pt*}

\subsection{A mod $p$ Langlands parametrization in families for $F=\bbQ_p$}

In this subsection we continue to assume that $F=\bbQ_p$ with $p\geq 5$. 

\begin{Pt*}
Recall the mod $p$ parametrization functor $P:\Mod(\cH^{(1)}_{\overline{\bbF}_p})\ra \SP_{\mathbf{\whG},0}$
from  \cite[7.3.6]{PS}. For $\zeta\in \cZ^{\vee}(k)$, let $\Mod_{\zeta}(\cH^{(1)}_{\overline{\bbF}_p})$ be the full subcategory of  $\Mod(\cH^{(1)}_{\overline{\bbF}_p})$ whose objets are the $\cH^{(1)}_{\overline{\bbF}_p}$-modules whose Satake parameter is supported on the closed subscheme $(V^{(1)}_{\mathbf{\whT},0}/W_0)_{\zeta}\subset V^{(1)}_{\mathbf{\whT},0}/W_0$.
A $\cH^{(1)}_{\overline{\bbF}_p}$-module $M$ lies in the category $\Mod_{\zeta}(\cH^{(1)}_{\overline{\bbF}_p})$ if and only if: $M$ is only supported in $\gamma$-components where 
$\gamma |_{\bbF_p^{\times}}=\zeta | _{\bbF_p^{\times}}$ and the operator $U^2$ acts on $M$ via the $\bbG_m$-part of 
$\zeta$.
Set 
$\SP_{\mathbf{\whG},0,\zeta}:=\QCoh((V^{(1)}_{\mathbf{\whT},0}/W_0)_{\zeta})$, the category of quasi-coherent modules on the $k$-scheme 
$(V^{(1)}_{\mathbf{\whT},0}/W_0)_{\zeta}$.
Then $P$ induces a \emph{mod $p$ $\zeta$-parametrization functor}
$$
\xymatrix{
P_{\zeta}:\Mod_{\zeta}(\cH^{(1)}_{\overline{\bbF}_p})\ar[r] & \SP_{\mathbf{\whG},0,\zeta}.
}
$$
 
For $\zeta\in \cZ^{\vee}(k)$, also recall the category $\LP_{\mathbf{\whG},0,\zeta}:=\QCoh(X_{\zeta})$ of mod $p$ Langlands parameters with determinant $\omega\zeta$ from \ref{modpLP}; it induces the functor 
$$
\xymatrix{
L_{\zeta*}:\SP_{\mathbf{\whG},0,\zeta}\ar[r] & \LP_{\mathbf{\whG},0,\zeta}
}
$$
\emph{push-forward along the $k$-morphism $L_{\zeta}:(V^{(1)}_{\mathbf{\whT},0}/W_0)_{\zeta}\ra X_{\zeta}$} from \ref{LLfamily}. 

Finally recall that for $\zeta\in \cZ^{\vee}(k)$, the functor of $I^{(1)}$-invariants $(\cdot)^{I^{(1)}}:\Mod^{\rm sm}(k[G])\ra \Mod(\cH^{(1)}_{\overline{\bbF}_p})$ induces a functor
$$
(\cdot)_{\zeta}^{I^{(1)}}:\Mod_{\zeta}^{\rm sm}(k[G])\ra \Mod_{\zeta}(\cH^{(1)}_{\overline{\bbF}_p}),
$$
by \ref{central_car_comp}.
\end{Pt*}

\begin{Def*}
Let $\zeta\in \cZ^{\vee}(k)$. The \emph{mod $p$ $\zeta$-Langlands parametrization functor} is the functor
$$
\LzPz:=L_{\zeta*}\circ P_{\zeta}:
$$
$$
\xymatrix{
\Mod_{\zeta}(\cH^{(1)}_{\overline{\bbF}_q})\ar[d] & \\
\LP_{\mathbf{\whG},0,\zeta}.
}
$$
Identifying $\zeta$ with a central character of $G$, the functor $\LzPz$ extends to the category $\Mod_{\zeta}^{\rm sm}(k[G])$ by precomposing with the functor $(\cdot)_{\zeta}^{I^{(1)}}:\Mod_{\zeta}^{\rm sm}(k[G])\ra \Mod_{\zeta}(\cH^{(1)}_{\overline{\bbF}_p})$:
$$
\LzPz\circ (\cdot)_{\zeta}^{I^{(1)}}:
$$
$$
\xymatrix{
\Mod_{\zeta}^{\rm sm}(k[G])\ar[d] & \\
\LP_{\mathbf{\whG},0,\zeta}.
}
$$
\end{Def*}

\begin{Th*}\label{paramfamily}
Suppose $F=\bbQ_p$ with $p\geq 5$. Fix a character $\zeta:  Z(G)=\bbQ_p^{\times} \ra k^\times$, corresponding to a point $(\zeta|_{\bbF_p^{\times}},\zeta(p^{-1}))\in\cZ^{\vee}(k)$ under the identification $\cZ(G)^{\vee}\cong\cZ^{\vee}(k)$ from \ref{ZGveeZvee}.

The mod $p$ $\zeta$-Langlands parametrization functor $\LzPz$ interpolates the Langlands parametrization of the blocks of the category 
$\Mod_{\zeta}^{\rm l adm}(k[G])$, cf. \ref{blocks} : for all $x\in X_{\zeta}(k)$ and for all $\pi\in b_{[\rho_x]}$,
$$
\LzPz(\pi^{I^{(1)}})=\left\{\begin{array}{ll} 
i_{x*}(\pi^{I^{(1)}}) & \textrm{if $x$ is not an exceptional point in the odd case} \\
i_{x*}(\pi^{I^{(1)}})^{\oplus 2} & otherwise 
\end{array}
\right.
\in \LP_{\mathbf{\whG},0,\zeta}
$$
where $i_x:\Spec(k)\ra X_{\zeta}$ is the $k$-point $x$.
\end{Th*}

\begin{proof}
By definition of a block of a category as a certain equivalence class of simple objects \cite{Pas13}, if $\pi\in b_{[\rho_x]}$ then in particular $\pi$ is simple. Then $\pi^{I^{(1)}}$ is simple too, and hence has a central character. Therefore $P_{\zeta}(\pi^{I^{(1)}})$ is the underlying $k$-vector space of $\pi^{I^{(1)}}$ supported at the $k$-point $v\in (V^{(1)}_{\mathbf{\whT},0}/W_0)_{\zeta}$ corresponding to its central character under the isomorphism $\sS^{(1)}_{\overline{\bbF}_p}$, which lies on some connected component $\gamma$. Suppose $\dim_k(\pi^{I^{(1)}})=2$. Then 
$\pi^{I^{(1)}}$ is isomorphic to the simple standard module of $\cH^{\gamma}_{\overline{\bbF}_p}$ with central character $v$, i.e. to 
$\ASph^{\gamma}(v)$, and hence $L_{\zeta}(v)=x$ by definition of the map of sets $L_{\zeta}(k)$. Suppose 
$\dim_k(\pi^{I^{(1)}})=1$. Then $\pi^{I^{(1)}}$ is one of the two antispherical characters of $\cH^{\gamma}_{\overline{\bbF}_p}$ whose restriction to the center $Z(\cH^{\gamma}_{\overline{\bbF}_p})$ is equal to $v$, i.e. it is one of the simple constituents of $(\ASph^{\gamma}(v))^{\sss}$, and hence again $L_{\zeta}(v)=x$ by definition of the map of sets $L_{\zeta}(k)$. Now if $x$ is not an exceptional point in an odd case, then $L_{\zeta}$ is an open immersion at $v$, and otherwise it has ramification index $2$ at $v$. The theorem follows.
\end{proof}

\vskip10pt 

\noindent {\small Cédric Pépin, LAGA, Université Paris 13, 99 avenue Jean-Baptiste Clément, 93 430 Villetaneuse, France \newline {\it E-mail address: \url{cpepin@math.univ-paris13.fr}} }

\vskip10pt

\noindent {\small Tobias Schmidt, IRMAR, Universit\'e de Rennes 1, Campus Beaulieu, 35042 Rennes, France \newline {\it E-mail address: \url{tobias.schmidt@univ-rennes1.fr}} }

\end{document}